\documentclass[10pt]{elsarticle}

\usepackage{amssymb,amsmath,amsthm,algorithmic,algorithm,booktabs,color}
\usepackage[mathscr]{eucal}

\newtheorem{theorem}{Theorem}[section]
\newtheorem{lemma}[theorem]{Lemma}

\newtheorem{proposition}[theorem]{Proposition}

\newtheorem{run}{Run}[section]

\newtheorem{remark}[theorem]{Remark}
\newtheorem{definition}[theorem]{Definition}

\DeclareMathOperator{\sgn}{sgn}

\journal{Advances in Computational Mathematics}

\begin{document}

\begin{frontmatter}



\title{Asymptotic Stability of POD based Model Predictive Control for a semilinear parabolic PDE}


\author{Alessandro Alla\corref{Thanks1}}
\ead{alessandro.alla@uni-hamburg.de}
\cortext[Thanks1]{This author wishes to acknowledge the support obtained by the ESF Grant no 4160.}
\address{Department of Mathematics, University of Hamburg, 20146 Hamburg, Germany}
\author{Stefan Volkwein\corref{Thanks2}}
\ead{stefan.volkwein@uni-konstanz.de}
\cortext[Thanks2]{This author gratefully acknowledges support by the DFG grant VO no 1658/2-1. S. Volkwein is the corresponding author.}
\address{Department of Mathematics and Statistics, University of Konstanz, 78457 Konstanz, Germany}

\begin{abstract}
In this article a stabilizing feedback control is computed for a semilinear parabolic partial differential equation utilizing a nonlinear model predictive (NMPC) method. In each level of the NMPC algorithm the finite time horizon open loop problem is solved by a reduced-order strategy based on proper orthogonal decomposition (POD). A stability analysis is derived for the combined POD-NMPC algorithm so that the lengths of the finite time horizons are chosen in order to ensure the asymptotic stability of the computed feedback controls. The proposed method is successfully tested by numerical examples. \end{abstract}


\end{frontmatter}

\section{Introduction}
\label{Section1}

In many control problems it is desired to design a stabilizing feedback control, but often the closed-loop solution can not be found analytically, even in the unconstrained case since it involves the solution of the corresponding Hamilton-Jacobi-Bellman equations; see, e.g., \cite{BCD97,Eva08} and \cite{KVX04}. But this approach requires the solution of a nonlinear hyperbolic partial differential equation with a high-dimensional spatial variable.

One approach to circumvent this problem is the repeated solution of open-loop optimal control problems. The first part of the resulting open-loop input signal is implemented and the whole process is repeated. Control approaches using this strategy are referred to as {\em model predictive control} (MPC), {\em moving horizon control} or {\em receding horizon control}. In general one distinguishes between linear and {\em nonlinear} MPC (NMPC). In linear MPC, linear models are used to predict the system dynamics and considers linear constraints on the states and inputs. Note that even if the system is linear, the closed loop dynamics are nonlinear due to the presence of constraints. NMPC refers to MPC schemes that are based on nonlinear models and/or consider a nonquadratic cost functional and general nonlinear constraints. Although linear MPC has become an increasingly popular control technique used in industry, in many applications linear models are not sufficient to describe the process dynamics adequately and nonlinear models must be applied. This inadequacy of linear models is one of the motivations for the increasing interest in nonlinear MPC; see. e.g., \cite {AFN04,FA01,GP11,RM09}. The prediction horizon plays a crucial role in MPC algorithms. For instance, the quasi infinite horizon NMPC allows an efficient formulation of NMPC while guaranteeing stability and the performances of the closed-loop as shown in \cite{AC98, FA02, IK02} under appropriate assumptions. For the purpose of our paper we will use a different approach since we will not deal with terminal constraints.

Since the computational complexity of NMPC schemes grows rapidly with the length of the optimization horizon, estimates for minimal stabilizing horizons are of particular interest to ensure stability while being computationally fast. Stability and suboptimality analysis for NMPC schemes without stabilizing constraints are studied in \cite[Chapter~6]{GP11}, where the authors give sufficient conditions ensuring asymptotic stability with minimal finite prediction horizon. Note that the stabilization of the problem and the computation of the minimal horizon involve the (relaxed) {\em dynamic programming principle} (DPP); see \cite{GPSW10, PRW11}. This approach allows estimates of the finite prediction horizon based on controllability properties of the dynamical system.

Since several optimization problems have to be solved in the NMPC method, it is reasonable to apply reduced-order methods to accelerate the NMPC algorithm. Here, we utilize proper orthogonal decomposition (POD) to derive reduced-order models for nonlinear dynamical systems; see, e.g., \cite{HLBR12,Sir87} and \cite{GV13}. The application of POD is justified by an a priori error analysis for the considered nonlinear dynamical system, where we combine techniques from \cite{KV01,KV02} and \cite{Sin13}. Let us refer to \cite{GU12}, where the authors also combine successfully an NMPC scheme with a POD reduced-order approach. However, no analysis is carried out ensuring the asymptotic stability of the proposed NMPC-POD scheme. Our contribution focusses on the stability analysis of the POD-NMPC algorithm without terminal constraints, where the dynamical system is a semilinear parabolic partial differential equation with an advection term. In particular, we study a minimal finite horizon for the reduced-order approximation such that it guarantees the asymptotic stability of the surrogate model. 
Our approach is motivated by the work \cite{AGW10}. The main difference here is that we have added an advection term in the dynamical system and utilize a POD suboptimal strategy to solve the open-loop problems. Since the minimal prediction horizon can be large, the numerical solution of the open-loop problems is very expensive within the NMPC algorithm. The application of the POD model  reduction reduces efficiently the computational cost by computing suboptimal solutions. But we involve this suboptimality in our stability analysis in order to ensure the asymptotic stability of our NMPC scheme.

The paper is organized in the following manner: In Section~\ref{Section2} we formulate our infinite horizon optimal control problem governed by a semilinear parabolic equation and bilateral control constraints. The NMPC algorithm is introduced in Section~\ref{Section3}. For the readers convenience, we recall the known results of the stability analysis. Further, the stability theory is applied to our underlying nonlinear semilinear equations and bilateral control constraints. In Section~\ref{Section4} we investigate the finite horizon open loop problem which has to be solved at each level of the NMPC algorithm. Moreover, we introduce the POD reduced-order approach and prove an a-priori error estimate for the semilinear parabolic equation. Finally, numerical examples are presented in Section~\ref{Section5}.

\section{Formulation of the control system}
\label{Section2}
\setcounter{equation}{0}
\setcounter{theorem}{0}
\setcounter{table}{0}
\setcounter{table}{0}

Let $\Omega=(0,1)\subset\mathbb R$ be the spatial domain. For the initial time $t_\circ\in\mathbb R^+_0=\{s\in\mathbb R\,|\,s\ge 0\}$ we define the space-time cylinder $Q=\Omega\times(t_\circ,\infty)$. By $H=L^2(\Omega)$ we denote the Lebesgue space of (equivalence classes of) functions which are (Lebesgue) measurable and square integrable. We endow $H$ by the standard inner product -- denoted by $\langle\cdot\,,\cdot\rangle_H$ -- and the associated induced norm $\|\varphi\|_H=\langle\varphi,\varphi\rangle_H^{1/2}$. Furthermore, $V=H^1_0(\Omega)\subset H$ stands for the Sobolev space
\[
V=\bigg\{\varphi\in H\,\Big|\,\int_\Omega\big|\varphi'(x)\big|^2\,\mathrm dx<\infty\text{ and }\varphi(0)=\varphi(1)=0\bigg\}.
\]
Recall that both $H$ and $V$ are Hilbert spaces. In $V$ we use the inner product
\[
{\langle\varphi,\phi\rangle}_V=\int_\Omega\varphi'(x)\phi'(x)\,\mathrm dx\quad\text{for }\varphi,\phi\in V
\]
and set $\|\varphi\|_V=\langle\varphi,\varphi\rangle_V^{1/2}$ for $\varphi\in V$. For more details on Lebesgue and Sobolev spaces we refer the reader to \cite{Eva08}, for instance. When the time $t$ is fixed for a given function $\varphi:Q\to\mathbb R$, the expression $\varphi(t)$ stands for a function $\varphi(\cdot\,,t)$ considered as a function in $\Omega$ only. Recall that the Hilbert space $L^2(Q)$ can be identified with the Bochner space $L^2(t_\circ,\infty;H)$. 

We consider the following control system governed by a semilinear parabolic partial differential equation: $y=y(x,t)$ solves the semilinear initial boundary value problem
\begin{subequations}
\label{PDESystem}
\begin{align}
\label{PDESystem-1}
y_t-\theta y_{xx}+y_x+\rho(y^3-y)&=u&&\text{in }Q,\\
\label{PDESystem-2}
y(0,\cdot)=y(1,\cdot)&=0 &&\text{in }(t_\circ,\infty),\\
\label{PDESystem-3}
y(t_\circ)&=y_\circ&&\text{in }\Omega.
\end{align}
\end{subequations}
In \eqref{PDESystem-1} it is assumed that the control $u=u(x,t)$ belongs to the set of admissible control inputs
\begin{equation}
\label{Uad(t0)}
\mathbb U_{ad}(t_\circ)=\big\{u\in\mathbb U(t_\circ)\,\big|\,u(x,t)\in U_{ad}\text{ for almost all (f.a.a.) } (x,t)\in Q\big\},
\end{equation}
where $\mathbb U(t_\circ)=L^2(t_\circ,\infty;H)$ and $U_{ad}=\{u\in\mathbb R\,|\,u_a\le u\le u_b\}$ with given $u_a\le0\le u_b$ . The parameters $\theta$ and $\rho$ satisfy
\[
(\theta,\rho)\in D_{ad}=\big\{(\tilde\theta,\tilde\rho)\in \mathbb R^2\,\big|\,\theta_a\le\tilde\theta\text{ and }\rho_a\le\tilde\rho\big\}
\]
with positive $\theta_a$ and $\rho_a$. Further, in \eqref{PDESystem-3} the initial condition $y_\circ=y_\circ(x)$ is supposed to belong to $H$.

A solution to \eqref{PDESystem} is interpreted in the weak sense as follows: for given $(t_\circ,y_\circ)\in \mathbb R^+_0\times H$ and $u\in\mathbb U_{ad}(t_\circ)$ we call $y$ a {\em weak solution} to \eqref{PDESystem} for fixed $(\theta,\rho)\in D_{ad}$ if $y(t)\in V$, $y_t(t)\in V'$ hold f.a.a. $t\ge t_\circ$ and $y$ satisfies $y(t_\circ)=y_\circ$ in $H$ as well as
\begin{equation}
\label{PDEWeakForm}
\frac{\mathrm d}{\mathrm dt}\,{\langle y(t),\varphi\rangle}_H+\int_\Omega\theta y_x(t)\varphi'+\big(y_x(t)+\rho(y^3(t)-y(t))\big)\varphi\,\mathrm dx=\int_\Omega u(t)\varphi\,\mathrm dx
\end{equation}
for all $\varphi\in V$ and f.a.a. $t>t_\circ$. Here, $y_t(t)$ stands for the distributional derivative with respect to the time variable satisfying \cite[p.~477]{DL00}
\[
\frac{\mathrm d}{\mathrm dt}\,{\langle y(t),\varphi\rangle}_H={\langle y_t(t),\varphi\rangle}_{V',V}\quad\text{for all }\varphi\in V.
\]
The following result is proved in \cite{CH98}, for instance.

\begin{proposition}
\label{Pro:ExUni}
For given $(t_\circ,y_\circ)\in \mathbb R^+_0\times H$ and $u\in \mathbb U_{ad}(t_\circ)$ there exists a unique weak solution $y=y_{[u,t_\circ,y_\circ]}$ to \eqref{PDESystem} for every $(\theta,\rho)\in D_{ad}$.
\end{proposition}

Let $(t_\circ,y_\circ)\in \mathbb R^+_0\times H$ be given. Due to Proposition~\ref{Pro:ExUni} we can define the quadratic cost functional:
\begin{equation}
\label{InfiniteHorizonCost}
\hat J(u;t_\circ,y_\circ):=\frac{1}{2}\int_{t_\circ}^\infty{\|y_{[u,t_\circ,y_\circ]}(t)-y_d\|}_H^2\,\mathrm dt+\frac{\lambda}{2}\int_{t_\circ}^\infty{\|u(t)\|}_H^2\,\mathrm dt
\end{equation}
for all $u\in \mathbb U(t_\circ)\supset \mathbb U_{ad}(t_\circ)$, where $y_{[u,t_\circ,y_\circ]}$ denotes the unique weak solution to \eqref{PDESystem}. We suppose that $y_d=y_d(x)$ is a given desired stationary state in $H$ (e.g., the equilibrium $y_d=0$) and that $\lambda>0$ denotes a fixed weighting parameter. Then we consider the nonlinear infinite horizon optimal control problem
\begin{equation}
\label{P}
\min \hat J(u;t_\circ,y_\circ)\quad\text{subject to (s.t.)}\quad u\in \mathbb U_{ad}(t_\circ).
\end{equation}

Suppose that the trajectory $y$ is measured at discrete time instances
\[
t_n=t_\circ+n\Delta t,\quad n\in\mathbb N,
\]
where the time step $\Delta t>0$ stands for the time step between two measurements. Thus, we want to select a control $u\in \mathbb U_{ad}(t)$ such that the associated trajectory $y_{[u,t_\circ,y_\circ]}$ follows a given desired state $y_d$ as good as possible. This problem is called a {\em tracking problem}, and, if $y_d=0$ holds, a {\em stabilization problem}.

Since our goal is to be able to react to the current deviation of the state $y$ at time $t=t_n$ from the given reference value $y_d$, we would like to have the control in {\em feedback form}, i.e., we want to determine a mapping $\mu:H\to\mathbb U_{ad}(t_\circ)$ with $u(t)=\mu(y(t))$ for $t\in[t_n,t_{n+1}]$.

\section{Nonlinear model predictive control}
\label{Section3}
\setcounter{theorem}{0}
\setcounter{equation}{0}
\setcounter{figure}{0}
\setcounter{table}{0}

We present an NMPC approach to compute a mapping $\mu$ which allows a representation of the control in feedback form. For more details we refer the reader to the monographs \cite{GP11,RM09}, for instance.
 
\subsection{The NMPC method}
\label{Section3.1}

To introduce the NMPC algorithm we write the weak form of our control system \eqref{PDESystem} as a parametrized nonlinear dynamical system. For $(\theta,\rho)\in D_{ad}$ let us introduce the $\theta$-and $\rho$-dependent nonlinear mapping $\mathcal F$ which maps the space $V\times H$ into the dual space $V'$ of $V$ as follows:
\[
\mathcal F(\varphi,v)=-\theta\varphi_{xx}+\varphi_x+\rho(\varphi^3-\varphi)-v\quad\text{for }(\varphi,v)\in V\times H.
\]
Then, we can express \eqref{PDEWeakForm} as the nonlinear dynamical system
\begin{equation}
\label{DynamicalSystem}
y'(t)=\mathcal F(y(t),u(t))\in V'\text{ for all }t>t_\circ,\quad y(t_\circ)=y_\circ\text{ in } H
\end{equation}
for given $(t_\circ,y_\circ)\in \mathbb R^+_0\times H$. The cost functional has been already introduced in \eqref{InfiniteHorizonCost}. Summarizing, we want to solve the following infinite horizon minimization problem
\begin{equation}
\tag{$\mathbf{P(}t_\circ\mathbf{)}$}
\label{P-DynSys}
\min \hat J(u;t_\circ,y_\circ)=\int_{t_\circ}^\infty \ell\big(y_{[u,t_\circ,y_\circ]}(t),u(t)\big)\,\mathrm dt\quad\text{s.t.}\quad u\in\mathbb U_{ad}(t_\circ),
\end{equation}
where we have defined the running quadratic cost as
\begin{equation}
\label{RunCost}
\ell(\varphi,v)=\frac{1}{2}\left({\|\varphi-y_d\|}^2_H+\lambda\,{\|v\|}_H^2\right)\quad\text{for }\varphi,v\in H.
\end{equation}
If we have determined a state feedback $\mu$ for \eqref{P-DynSys}, the control $u(t)=\mu(y(t))$ allows a closed loop representation for $t\in[t_\circ,\infty)$. Then, for a given initial condition $y_0\in H$ we set $t_\circ=0$, $y_\circ=y_0$ in \eqref{DynamicalSystem} and insert $\mu$ to obtain the closed-loop form
\begin{equation}
\label{DynSysClosLoop}
\begin{aligned}
&y'(t)=\mathcal F(y(t),\mu(y(t)))&&\text{in } V' \text{ for }t\in(t_\circ,\infty),\\
&y(t_\circ)=y_\circ&&\text{in } H.
\end{aligned}
\end{equation}
Note that the infinite horizon problem may be very hard to solve due to the dimensionality of the problem. On the other hand it guarantees the stabilization of the problem which is very important for certain applications. In an NMPC algorithm a state feedback law is computed for \eqref{P-DynSys} by solving a sequence of finite time horizon problems.

To formulate the NMPC algorithm we introduce the finite horizon quadratic cost functional as follows: for $(t_\circ,y_\circ)\in \mathbb R^+_0\times H$ and $u\in \mathbb U_{ad}^N(t_\circ)$ we set
\[
\hat J^N(u;t_\circ,y_\circ)=\int_{t_\circ}^{t_\circ^N} \ell\big(y_{[u,t_\circ,y_\circ]}(t),u(t)\big)\,\mathrm dt,
\]
where $N$ is a natural number, $t_\circ^N=t_\circ+N\Delta t$ is the final time and $N\Delta t$ denotes the length of the time horizon for the chosen time step $\Delta t>0$. Further, we introduce the Hilbert space $\mathbb U^N(t_\circ)=L^2(t_\circ,t_\circ^N;H)$ and the set of admissible controls
\[
\mathbb U_{ad}^N(t_\circ)=\big\{u\in \mathbb U^N(t_\circ)\,\big|\,u(x,t)\in U_{ad}\text{ f.a.a. } (x,t)\in Q^N\big\}
\]
with $Q^N=\Omega\times(t_\circ,t_\circ^N)\subset Q$; compare \eqref{Uad(t0)}. In Algorithm~\ref{Alg:NMPC} the method is presented.
\begin{algorithm}
\caption{(NMPC algorithm)}
\label{Alg:NMPC}
\begin{algorithmic}[1]
\REQUIRE time step $\Delta t>0$, finite horizon $N\in\mathbb N$, weighting parameter $\lambda>0$.
\FOR{$n=0,1,2,\ldots$}
\STATE Measure the state $y(t_n)\in V$ of the system at $t_n=n\Delta t$.
\STATE Set $t_\circ=t_n=n\Delta t$, $y_\circ=y(t_n)$ and compute a global solution to 
\begin{equation}
\tag{$\mathbf{P^{\boldsymbol N}(}t_\circ\mathbf{)}$}
\label{P-DynSys-FinHor}
\min\hat J^N(u;t_\circ,y_\circ)\quad\text{s.t.}\quad u\in\mathbb U_{ad}^N(t_\circ).
\end{equation}
We denote the obtained optimal control by $\bar u^N$.
\STATE Define the NMPC feedback value $\mu^N(t;t_\circ,y_\circ)=\bar u^N(t)$, $t\in(t_\circ,t_\circ+\Delta t]$ and use this control to compute the associated state $y=y_{[\mu^N(\cdot),t_\circ,y_\circ]}$ by solving \eqref{DynamicalSystem} on $[t_\circ,t_\circ+\Delta t]$.
\ENDFOR
\end{algorithmic}
\end{algorithm}

\noindent
We store the optimal control on the first subinterval $[t_\circ,t_\circ+\Delta t]=[0,\Delta t]$ and the associated optimal trajectory. Then, we initialize a new finite horizon optimal control problem whose initial condition is given by the optimal trajectory $\bar y(t)=y_{[\mu^N(\cdot),t_\circ,y_\circ]}(t)$ at $t=t_\circ+\Delta t$ using the optimal control $\mu^N(t;t_\circ,y_\circ)=\bar u^N(t)$ for $t\in(t_\circ,t_\circ+\Delta t]$ . We iterate this process by setting $t_\circ=t_\circ+\Delta t$. Of course, the larger the horizon, the better the approximation one can have, but we would like to have the minimal horizon which can guarantee stability \cite{GPSW10}. Note that \eqref{P-DynSys-FinHor} is an open loop problem on a finite time horizon $[t_\circ,t_\circ+N\Delta t]$ which will be studied in Section~\ref{Section4}. 

\subsection{Dynamic programming principle (DPP) and asymptotic stability}\label{Section3.2}

For the reader's convenience we now recall the essential theoretical results from dynamic programming and stability analysis. Let us first introduce the so called {\em value function} $v$ defined as follows for an infinite horizon optimal control problem:
\[
v(t_\circ,y_\circ):=\inf_{u\in\mathbb U_{ad}(t_\circ)}\hat J(u;t_\circ,y_\circ)\quad\text{for }(t_\circ,y_\circ)\in \mathbb R_0^+\times H.
\]
Let $N\in\mathbb N$ be chosen. The DDP states that the value function $v$ satisfies for any $k\in\{1,\ldots,N\}$ with $t_\circ^k=t_k+k\Delta t$:
\begin{align*}
&v(t_\circ,y_\circ)\\
&=\inf_{u\in\mathbb U_{ad}^k(t_\circ)} \left\{\int_{t_\circ}^{t_\circ^k}\ell\big(y_{[u,t_\circ,y_\circ]}(t),u(t)\big)\,\mathrm dt + v\big(t_\circ+k\Delta t,y_{[u,t_\circ,y_\circ]}(t_\circ+k\Delta t)\big)\right\}
\end{align*}
which holds under very general conditions on the data; see, e.g., \cite{BCD97} for more details. The value function for the finite horizon problem \eqref{P-DynSys-FinHor} is of the following form:
\[
v^N(t_\circ,y_\circ)=\inf_{u\in\mathbb U_{ad}^N(t_\circ)}\hat J^N(u;t_\circ,y_\circ)\quad\text{for }(t_\circ,y_\circ)\in \mathbb R^+_0\times H.
\]
The value function $v^N$ satisfies the DPP for the finite horizon problem for $t_\circ+k\Delta t$, $0<k<N$:
\begin{align*}
&v^N(t_\circ,y_\circ)\\
&\quad=\inf_{u\in\mathbb U_{ad}^k(t_\circ)}\left\{\int_{t_\circ}^{t_\circ+k\Delta t} \ell\big(y_{[u,t_\circ,y_\circ]}(t),u(t)\big)\,\mathrm dt+ v^{N-k}\big(y_{[u,t_\circ,y_\circ]}(t_\circ+k\Delta t)\big)\right\}.
\end{align*}
Nonlinear stability properties can be expressed by comparison functions which we recall here for the readers convenience \cite[Definition~2.13]{GP11}.

\begin{definition}
\label{De:2.13}
We define the following classes of comparison functions:
\begin{align*}
\mathcal K &= \big\{\beta:\mathbb R^+_0\to\mathbb R^+_0\,\big|\,\beta\text{ is continuous, strictly increasing and } \beta(0)=0\big\},\\
\mathcal K_\infty &= \big\{\beta:\mathbb R^+_0\to\mathbb R^+_0\,\big|\,\beta\in\mathcal K,\,\beta \text{ is unbounded}\big\},\\
\mathcal L &= \left\{\beta:\mathbb R^+_0\to\mathbb R^+_0\,\big|\,\beta\text{ is continuous, strictly decreasing, } \lim\limits_{t\to\infty}
\beta(t)=0\right\},\\
\mathcal K\mathcal L &= \big\{\beta:\mathbb R^+_0\times\mathbb R^+_0\to\mathbb R^+_0\,
\big|\,\beta\text{ is continuous, }\beta(\cdot\,,t)\in\mathcal K,\,\beta(r,
\cdot)\in\mathcal L\big\}.
\end{align*}
\end{definition}

Utilizing a comparison function $\beta\in\mathcal K\mathcal L$ we introduce the concept of asymptotic stability; see, e.g. \cite[Definition~2.14]{GP11}.

\begin{definition}
\label{De:2.14}
Let  $y_{[\mu(\cdot),t_\circ,y_\circ]}$ be the solution to \eqref{DynSysClosLoop} and $y_*\in H$ an equilibrium for \eqref{DynSysClosLoop}, i.e., we have $\mathcal F(y_*,\mu(y_*))=0$. Then, $y_*$ is said to be {\em locally asymptotically stable} if there exist a constant $\eta>0$ and a function $\beta\in\mathcal K\mathcal L$ such that the estimate
\begin{equation*}
{\|y_{[\mu(\cdot),t_\circ,y_\circ]}(t)-y_*\|}_H \le \beta\big({\|y_\circ-y_*\|}_H,t)
\end{equation*}
holds for all $y_\circ\in H$ satisfying $\|y_\circ-y_*\|_H<\eta$ and all $t\ge t_\circ$.
\end{definition}

Let us recall the main result about asymptotic stability via DPP; see \cite{GPSW10}.

\begin{proposition}
Let $N\in\mathbb N$ be chosen and the feedback mapping $\mu^N$ be computed by Algorithm~{\rm\ref{Alg:NMPC}}. Assume that there exists an $\alpha^N\in (0,1]$ such that for all $(t_\circ,y_\circ)\in \mathbb R^+_0\times H$ the {\em relaxed DPP}
\begin{equation}
\label{rel:dpp}
v^N(t_\circ,y_\circ)\ge v^N\big(t_\circ+\Delta t,y_{[\mu^N(\cdot),t_\circ,y_\circ]}(t_\circ+\Delta t)\big)+\alpha^N \ell\big(y_\circ,\mu^N(y_\circ))\big)
\end{equation}
holds. Then we have for all $(t_\circ,y_\circ)\in \mathbb R^+_0\times H$:
\begin{equation}
\label{rel:est}
\alpha^N v(t_\circ,y_\circ)\le\alpha^N \hat J(\mu^N(y_{[\mu^N(\cdot),t_\circ,y_\circ]});t_\circ,y_\circ)\leq v^N(t_\circ,y_\circ)\le v(t_\circ,y_\circ),
\end{equation}
where $y_{[\mu^N(\cdot),t_\circ,y_\circ]}$ solves the closed-loop dynamics \eqref{DynSysClosLoop} with $\mu=\mu^N$. If, in addition, there exists an equilibrium $y_*\in H$ and $\alpha_1,\alpha_2\in\mathcal K_\infty$ satisfying
\begin{subequations}
\label{growth}
\begin{align}
\label{growth-1}
&\ell_*(y_\circ)=\min_{u\in U_{ad}} \ell(y_\circ,u)\ge \alpha_1\big({\|y_\circ-y_*\|}_H\big),\\
\label{growth-2}
&\alpha_2\big({\|y_\circ-y_*\|}_H\big)\ge v^N(t_\circ,y_\circ)
\end{align}
\end{subequations}
hold for all $(t_\circ,y_\circ)\in \mathbb R^+_0\times H$, then $y_*$ is a globally {\em asymptotically stable equilibrium} for \eqref{DynSysClosLoop} with the feedback map $\mu=\mu^N$ and value function $v^N$.
\end{proposition}

\begin{remark}
\label{rmk_imp}
\rm
\begin{enumerate}
\item [1)] Our running cost $\ell$ defined in \eqref{RunCost} satisfies condition \eqref{growth-1} for the choice $y_d=y_*$. Further, \eqref{growth-2} follows from the finite horizon quadratic cost functional $\hat J^N$, the definition of the value function $v^N$ and our a-priori analysis presented in Lemma~\ref{lem:stab} below. Therefore, we only have to check the relaxed DPP \eqref{rel:dpp}.
\item [2)] It is proved in \cite {GPSW10} that $\lim\limits_{N\to\infty} \alpha^N=1$. Hence, we would like to find $\alpha^N$ close to one to have the best approximation of $v$ in terms of $v^N$. On the other hand, a large $N$ implies that the numerical solution of \eqref{P-DynSys-FinHor} is much more involved. We will discuss the numerical computation of $\alpha^N$ next.
\item [3)] By \eqref{rel:est} we obtain the suboptimality estimate
\[
\hat J(\mu^N\big(y_{[\mu^N(\cdot),t_\circ,y_\circ]});t_\circ,y_\circ\big)\le \frac{v^N(t_\circ,y_\circ)}{\alpha^N}\le \frac{v(t_\circ,y_\circ)}{\alpha^N};
\]
compare \cite[Section~4.3]{GP11}.\hfill$\Diamond$
\end{enumerate}
\end{remark}

\noindent
In order to estimate $\alpha^N$ in the relaxed DPP we require the exponential controllability property for the system.

\begin{definition}
\label{De:ExpSta}
System \eqref{DynamicalSystem} is called {\em exponentially controllable} with respect to the running cost $\ell$ if for each $(t_\circ,y_\circ)\in \mathbb R^+_0\times H$ there exist two real constants $C>0$, $\sigma\in[0,1)$ and an admissible control $u\in\mathbb U_{ad}(t_\circ)$ such that:%
\begin{equation}
\label{exp:con}
\ell(y_{[u,t_\circ,y_\circ]}(t),u(t))\leq C\sigma^{t-t_\circ}\ell_*(y_\circ)\quad\text{f.a.a. } t\ge t_\circ.
\end{equation}
\end{definition}

We present an a-priori estimate for the uncontrolled solution to \eqref{DynamicalSystem}, i.e., the solution for $u=0$. For a proof we refer to the \ref{App-Lemma3.6}. Recall that $V$ is continuously (even compactly) embedded into $H$. Due to the Poincar\'e inequality \cite{Eva08} there exists a constant $C_V>0$ such that
\begin{equation}
\label{Eq-ocp.embed}
{\|\varphi\|}_H\le C_V\,{\|\varphi\|}_V\quad\text{for all }\varphi\in V.
\end{equation}

\begin{lemma}
\label{lem:stab}
Let $(t_\circ,y_\circ)\in \mathbb R_0^+\times H$ and $u=-Ky\in\mathbb U_{ad}(t_\circ)$ with an appropriate real constant $K>0$. Then, the solution $y=y_{[u,t_\circ,y_\circ]}$ to \eqref{DynamicalSystem} satisfies the a-priori estimate
\begin{equation}
\label{AprioriEst}
{\|y(t)\|}_H\le e^{-\gamma(K)(t-t_\circ)}\,{\|y_\circ\|}_H\quad\text{f.a.a. }t\ge t_\circ
\end{equation}
with $\gamma(K)=\gamma(K;\theta,\rho)=K+\theta/C_V-\rho$.
\end{lemma}

\begin{remark}
\label{Re:AprioriEst}
\rm
\begin{enumerate}
\item [1)] Let $K=0$ hold. Then, for $\theta>\rho C_V$ we have $\gamma>0$. Then, \eqref{AprioriEst} implies that $\|y(t)\|_H<\|y_\circ\|_H$ for any $t>t_\circ$. Moreover, the origin $y_\circ=0$ is unstable for $\gamma<0$; see\cite[Example~6.27]{GP11}.
\item [2)] If $K>\rho-\theta/C_V$ holds, $\|y(t)\|_H$ tends to zero for $t\to\infty$.\hfill$\Diamond$
\end{enumerate}
\end{remark}

Let us choose $y_d=0$. Suppose that we have a particular class of state feedback controls of the form $u(x,t)=-Ky(x,t)$ with a positive constant $K$; see \cite{AGW10}. This assumption helps us to derive the exponential controllability in terms of the running cost $\ell$ and to compute a minimal finite time prediction horizon $N\Delta t$ ensuring asymptotic stability. Combining \eqref{AprioriEst} with the desired exponential controllability \eqref{exp:con} and using $y_d=0$ we obtain for all $t\ge t_\circ$ \cite{AGW10}:
\begin{equation}\label{stab:prp}
\begin{aligned}
\ell(y(t),u(t))&=\frac{1}{2}\,\big({\|y(t)\|}_H^2+\lambda\,{\|u(t)\|}_H^2\big)=\frac{1}{2}\,(1+\lambda K^2)\,{\|y(t)\|}_H^2\\
&\le\frac{1}{2}\,C(K)e^{-2\gamma(K)(t-t_\circ)}\,{\|y_\circ\|}_H^2=C(K)\sigma(K)^{t-t_\circ}\,\ell_*(y_\circ)
\end{aligned}
\end{equation}
f.a.a. $t\ge t_\circ$ and for every $(t_\circ,y_\circ)\in \mathbb R^+_0\times H$, where
\begin{equation}
\label{C(K)}
C(K)=(1+\lambda K^2),\quad \sigma(K)=e^{-2\gamma(K)}.
\end{equation}

In the following theorem we provide an explicit formula for the scalar $\alpha^N$ in \eqref{rel:dpp}. A complete discussion is given in \cite{GPSW10}.

\begin{theorem}
\label{Th:Alpha}
Assume that the system \eqref{DynamicalSystem} and $\ell$ statisfy the controllability condition \eqref{exp:con}. Let the finite prediction horizon $N\Delta t$ be given with $N\in\mathbb N$ and $\Delta t>0$. Then the parameter $\alpha^N$ depends on $K$ and is given by:
\begin{equation}
\label{alpha}
\alpha^N(K)=1-\frac{\big(\eta_N(K)-1\big)\prod_{i=2}^N\big(\eta_i(K)-1\big)}{\prod_{i=2}^N\eta_i(K)- \prod_{i=2}^N\big(\eta_i(K)-1\big)}
\end{equation}
where $\eta_i(K)=C(1-\sigma^i)/(1-\sigma)$ and the constants $C=C(K)$, $\sigma=\sigma(K)$ are given by \eqref{C(K)}.
\end{theorem}

\begin{remark}
\label{Rem:Kest}
\rm
\begin{enumerate}
\item [1)] Theorem~\ref{Th:Alpha} suggests how we can compute a minimal horizon $N$ which ensures asympotic stability; see \cite{AG13}. Due to \eqref{C(K)} we fix a small finite horizon $N\in\mathbb N$ compute a (global) solution $\bar K$ to
\begin{equation}
\label{alphamax}
\max \alpha^N(K)\quad\text{s.t.}\quad \gamma(K) \ge \varepsilon
\end{equation}
with $0<\varepsilon\ll1$ and $\eta_i(K)$ from Theorem~\ref{Th:Alpha}. If the optimal value $\alpha^N(\bar K)$ is greater than zero, the finite horizon guarantees asymptotic stability. If $\alpha^N(\bar K)<0$ holds, we enlarge $N$ and solve \eqref{alphamax} again.
\item [2)] Since we suppose that $u\in\mathbb U^N_{ad}(t_\circ)$, we have to guarantee the bilateral control constraints
\begin{equation}
\label{Eq:K-uab}
u_a\leq-Ky(x,t)\leq u_b\quad\text{f.a.a. }(x,t)\in Q^N
\end{equation}
with $u_a\le0\le u_b$. This leads to additional constraints for $K$ in \eqref{alphamax}. Since we determine $K$ in such a way that $\gamma(K)>0$ is satisfied, we derive from \eqref{AprioriEst} that
\[
{\|y(t)\|}_H\le{\|y_\circ\|}_H\quad\text{f.a.a. }t\ge t_\circ.
\]
Let us suppose that we have $y_\circ\neq0$ and $\|y(t)\|_{C(\overline\Omega)}\le \|y_\circ\|_{C(\overline\Omega)}$ f.a.a. $t\ge t_\circ$. Then, we define
\begin{equation}
\label{Eq:K-y0}
y_{\circ a}=\min_{x\in\overline \Omega} y_\circ(x),\quad y_{\circ b}=\max_{x\in\overline \Omega} y_\circ(x).
\end{equation}
Then, $K$ has to satisfy $\gamma(K) \ge \varepsilon$ and the restrictions shown in Table~\ref{K:cond1}.
\begin{table}
\begin{center}
\begin{tabular}{c|cc}
\toprule
$K$ & $y_{\circ a}<0$ & $y_{\circ a}\ge0$ \\\toprule
$y_{\circ b}=0$ & no constraints & not considered\\\midrule
$y_{\circ b}<0$ & $K\le u_b/|y_{\circ b}|$ & impossible\\\midrule
$y_{\circ b}>0$ & $K\le\min\big\{|u_a|/y_{\circ b},u_b/|y_{\circ a}|\big\} $& $K\le|u_a|/y_{\circ b}$\\
\bottomrule
\end{tabular}
\end{center}
\caption{Constraints for the feedback factor $K$ in $u(x,t)=-Ky(x,t)$ considering the bilateral control constraints \eqref{Eq:K-uab} and the initial condition \eqref{Eq:K-y0}.}
\label{K:cond1}
\end{table}
Summarizing, $K$ has always an upper bound due to the constraints $u_a$, $u_b$ and a lower bound due to the stabilization related to $\gamma(K)>0$. \hfill$\Diamond$
\end{enumerate}
\end{remark}

\section{The finite horizon problem \eqref{P-DynSys-FinHor}}
\label{Section4}
\setcounter{equation}{0}
\setcounter{theorem}{0}
\setcounter{figure}{0}
\setcounter{table}{0}

In this section we discuss \eqref{P-DynSys-FinHor}, which has to be solved at each level of Algorithm~\ref{Alg:NMPC}.

\subsection{The open loop problem}
\label{Section4.1}

Recall that we have introduced the final time $t_\circ^N=t_\circ+N\Delta t$ and the control space $\mathbb U^N(t_\circ)=L^2(t_\circ,t_\circ^N;H)$. The space $\mathbb Y^N(t_\circ)=W(t_\circ,t_\circ^N)$ is given by
\[
W(t_\circ,t_\circ^N)=\big\{\varphi\in L^2(t_\circ,t_\circ^N;V)\,\big|\,\varphi_t\in L^2(t_\circ,t_\circ^N;V')\big\},
\]
which is a Hilbert space endowed with the common inner product \cite[pp.~472-479]{DL00}. We define the Hilbert space $\mathbb X^N(t_\circ)=\mathbb Y^N(t_\circ)\times\mathbb U^N(t_\circ)$ endowed with the standard product topology. Moreover, we introduce the Hilbert space $\mathbb Z^N(t_\circ)=\mathbb Z^N_1(t_\circ)\times H$ with $\mathbb Z^N_1(t_\circ)=L^2(t_\circ,t_\circ^N;V)$ and the nonlinear operator $e=(e_1,e_2):\mathbb X^N(t_\circ)\to \mathbb Z^N(t_\circ)'$ by
\begin{align*}
&{\langle e_1(x),\varphi\rangle}_{\mathbb Z^N_1(t_\circ)',\mathbb Z^N_1(t_\circ)}=\int_{t_\circ}^{t_\circ^N}{\langle y_t(t),\varphi(t)\rangle}_{V',V}\,\mathrm dt\\
&\quad+\int_{t_\circ}^{t_\circ^N}\int_\Omega \theta y_x(t)\varphi(x)+\Big(y_x(t)+\rho\big(y(t)^3-y(t)\big)-u(t)\Big)\varphi(t)\,\mathrm dx\mathrm dt,\\
&{\langle e_2(x),\phi\rangle}_H={\langle y(t_\circ)-y_\circ,\phi\rangle}_H
\end{align*}
for $x=(y,u)\in\mathbb X^N(t_\circ)$, $(\varphi,\phi)\in \mathbb Z^N(t_\circ)$, where we identify the dual $\mathbb Z^N(t_\circ)'$ of $\mathbb Z^N(t_\circ)$ with $L^2(t_\circ,t_\circ^N;V')\times H$ and $\langle\cdot\,,\cdot\rangle_{\mathbb Z^N_1(t_\circ)',\mathbb Z^N_1(t_\circ)}$ denotes the dual pairing between $\mathbb Z^N_1(t_\circ)'$ and $\mathbb Z^N_1(t_\circ)$. Then, for given $u\in\mathbb U^N(t_\circ)$ the weak formulation for \eqref{PDEWeakForm} can be expressed as the operator equation $e(x)=0$ in $\mathbb Z^N(t_\circ)'$. Further, we can write \eqref{P-DynSys-FinHor} as a constrained infinite dimensional minimization problem
\begin{equation}
\label{PN-OptPro}
\min J(x)=\int_{t_\circ}^{t_\circ^N}\ell(y(t),u(t))\,\mathrm dt\quad\text{s.t.}\quad x\in\mathbb F_{ad}^N(t_\circ)
\end{equation}
with the feasible set
\[
\mathbb F_{ad}^N(t_\circ)=\big\{x=(y,u)\in \mathbb X^N(t_\circ)\,\big|\,e(x)=0\text{ in }\mathbb Z^N(t_\circ)'\text{ and }u\in\mathbb U^N_{ad}(t_\circ)\big\}.
\]
For given fixed control $u\in\mathbb U_{ad}^N(t_\circ)$ we consider the state equation $e(y,u)=0\in \mathbb Z^N(t_\circ)'$, i.e., $y$ satisfies
\begin{equation}
\label{StateEquation}
\begin{aligned}
&\frac{\mathrm d}{\mathrm dt}\,{\langle y(t),\varphi\rangle}_H+\int_\Omega\theta y_x(t)\varphi'+\big(y_x(t)+\rho(y(t)^3-y(t))\big)\varphi\,\mathrm dx\\
&\hspace{50mm}=\int_\Omega u(t)\varphi\,\mathrm dx\quad\text{f.a.a. }t\in(t_\circ,t_\circ^N],\\
&{\langle y(t_\circ),\varphi\rangle}_H={\langle y_\circ,\varphi\rangle}_H
\end{aligned}
\end{equation}
for all $\varphi\in V$. The following result is proved in \cite[Theorem~5.5]{Tro10}.

\begin{proposition}
\label{Pro:StateEq}
For given $(t_\circ,y_\circ)\in \mathbb R^+_0\times H$ and $u\in \mathbb U_{ad}^N(t_\circ)$ there exists a unique weak solution $y\in \mathbb Y^N(t_\circ)$ to \eqref{StateEquation} for every $(\theta,\rho)\in D_{ad}$. If, in addition, $y_\circ$ is essentially bounded in $\Omega$, i.e., $y_\circ\in L^\infty(\Omega)$ holds, we have $y\in L^\infty(Q^N)$ satisfying
\begin{equation}
\label{APrioriState}
{\|y\|}_{\mathbb Y^N(t_\circ)}+{\|y\|}_{L^\infty(Q^N)}\le C\big({\|u\|}_{\mathbb U^N(t_\circ)}+{\|y_\circ\|}_{L^\infty(\Omega)}\big)
\end{equation}
for a $C>0$, which is independent of $u$ and $y_\circ$.
\end{proposition}

Utilizing \eqref{APrioriState} it can be shown that \eqref{PN-OptPro} possesses at least one (local) optimal solution which we denote by $\bar x^N=(\bar y^N,\bar u^N)\in \mathbb F_{ad}^N(t_\circ)$; see \cite[Chapter~5]{Tro10}. For the numerical computation of $\bar x^N$ we turn to first-order necessary optimality conditions for \eqref{PN-OptPro}. To ensure the existence of a unique Lagrange multiplier we investigate the surjectivity of the linearization $e'(\bar x^N):\mathbb X^N(t_\circ)\to \mathbb Z^N(t_\circ)'$ of the operator $e$ at a given point $\bar x^N=(\bar y^N,\bar u^N)\in \mathbb X^N(t_\circ)$. Note that the Fr\'echet derivative $e'(\bar x^N)=(e_1'(\bar x^N),e_2'(\bar x^N))$ of $e$ at $\bar x^N$ is given by
\begin{align*}
&{\langle e_1'(\bar x^N)x,\varphi\rangle}_{\mathbb Z_1^N(t_\circ)',\mathbb Z_1^N(t_\circ)}=\int_{t_\circ}^{t_\circ^N}{\langle y_t(t),\varphi(t)\rangle}_{V',V}\,\mathrm dt\\
&\qquad+\int_{t_\circ}^{t_\circ^N}\int_\Omega \theta y_x(t)\varphi(x)+\Big(y_x(t)+\rho\big(3\bar y^N(t)^2-1\big)y(t)-u(t)\Big)\varphi(t)\,\mathrm dx\mathrm dt,\\
&{\langle e_2'(\bar x^N)x,\phi\rangle}_H={\langle y(t_\circ),\phi\rangle}_H
\end{align*}
for $x=(y,u)\in\mathbb X^N(t_\circ)$, $(\varphi,\phi)\in \mathbb Z^N(t_\circ)$. Now, the operator $e'(\bar x^N)$ is surjective if and only if for an arbitrary $F=(F_1, F_2)\in\mathbb Z^N(t_\circ)'$ there exists a pair $x=(y,u)\in \mathbb X^N(t_\circ)$ satisfying $e'(\bar x^N)=F$ in $\mathbb Z^N(t_\circ)'$ which is equivalent with the fact that there exist a $u\in \mathbb U^N(t_\circ)$ and a $y\in\mathbb Y^N(t_\circ)$ solving the linear parabolic problem
\begin{equation}
\label{e'(x)}
y_t-\theta y_{xx}+y_x+\rho (3\bar{y}^2-1)y=F_1\text{ in }\mathbb Z_1^N(t_\circ)',\quad y(t_\circ)=F_2\text{ in }H.
\end{equation}
Utilizing standard arguments \cite{DL00} it follows that there exists for any $u\in \mathbb U^N(t_\circ)$ a unique $y\in \mathbb Y^N(t_\circ)$ solving \eqref{e'(x)}. Thus, $e'(\bar x^N)$ is a surjective operator and the local solution $\bar x^N$ to \eqref{PN-OptPro} can be characterized by first-order optimality conditions. We introduce the Lagrangian by
\[
L(x,p,p_\circ)=J(x)+{\langle e(x), (p,p_\circ)\rangle}_{\mathbb Z^N(t_\circ)',\mathbb Z^N(t_\circ)}\]
for $x\in \mathbb X^N(t_\circ)$ and $(p,p_\circ)\in \mathbb Z^N(t_\circ)$. Then, there exists a unique associated Lagrange multiplier pair $(\bar p^N,\bar p_\circ)$ to \eqref{PN-OptPro} satisfying the optimality system
\begin{align*}
&\nabla_y L(\bar x^N,\bar p^N,\bar p_\circ^N)y=0&&\forall y\in \mathbb Y^N(t_\circ)&&\hspace{-3mm}\text{\em (adjoint equation)}\\
&\nabla_u L(\bar x^N,\bar p^N,\bar p_\circ^N)(u-\bar u^N) \ge 0&&\forall u\in\mathbb U_{ad}^N(t_\circ)&&\hspace{-3mm}\text{\em (variational inequality)},\\
&{\langle e(\bar x^N),(p,p_\circ)\rangle}_{\mathbb Z^N(t_\circ)',\mathbb Z^N(t_\circ)}=0&&\forall(\bar{p},\bar{p}_0)\in \mathbb Z^N(t_\circ)&&\hspace{-3mm}\text{\em (state equation)}.
\end{align*}
It follows from variational arguments that the strong formulation for the adjoint equation is of the form
\begin{equation}
\label{DualEq}
\begin{aligned}
-\bar p^N_t-\theta \bar p^N_{xx}-\bar p^N_x-\rho\big(1-3(\bar y^N)^2\big)\bar p^N&=y_d-\bar y^N&&\text{in }Q^N,\\
\bar p^N(0,\cdot)=\bar p^N(1,\cdot)&=0&&\text{in }(t_\circ,t_\circ^N),\\
\bar p^N(t_\circ^N)&=0&&\text{in }\Omega.
\end{aligned}
\end{equation}
Moreover, we have $\bar p_\circ^N=\bar p^N(t_\circ)$. The variational inequality base the form
\begin{equation}
\label{VarIneq}
\int_{t_\circ}^{t_\circ^N}\int_\Omega(\lambda \bar u^N-\bar p^N)(u-\bar u^N)\,\mathrm dx\mathrm dt\ge 0\quad \text{for all }u\in\mathbb U^N_{ad}(t_\circ).
\end{equation}

Using the techniques as in \cite[Proposition~2.12]{Vol03} one can prove that second-order sufficient optimality conditions can be ensured provided the residuum $\|\bar y^N-y_d\|_{L^2(t_\circ,t^N_\circ;H)}$ is sufficiently small.

\subsection{POD reduced order model for open-loop problem}
\label{Section4.2}

To solve \eqref{PN-OptPro} we apply a reduced-order discretization based on proper orthogonal decomposition (POD); see \cite{GV13}. In this subsection we briefly introduce the POD method, present an a-priori error estimate for the POD solution to the state equation $e(x)=0\in \mathbb Z^N(t_\circ)'$ and formulate the POD Galerkin approach for \eqref{PN-OptPro}.

\subsubsection{The POD method for dynamical systems}
\label{Section4.2.1}

By $X$ we denote either the function space $H$ or $V$. Then, for $\wp\in\mathbb N$ let the so-called {\em snapshots} or {\em trajectories} $y^k(t)\in X$ be given f.a.a. $t\in [t_\circ,t_\circ^N]$ and for $1\le k\le\wp$. At least one of the trajectories $y^k$ is assumed to be nonzero. Then we introduce the linear subspace
\begin{equation}
\label{Eq:SnapshotSpace}
\mathscr V=\mathrm{span}\,\Big\{y^k(t)\,|\,t\in [t_\circ,t_\circ^N]\text{ a.e. and } 1\le k\le \wp\Big\}\subset X
\end{equation}
with dimension $d\ge 1$. We call the set $\mathscr V$ {\em snapshot subspace}. The method of POD consists in choosing a complete orthonormal basis in $X$ such that for every $\mathfrak l \le d$ the mean square error between $y^k(t)$ and their corresponding $\mathfrak l$-th partial Fourier sum is minimized on average:
\begin{equation}
\tag{$\mathbf P^\mathfrak l$}
\label{Eq:PODMin}
\left\{
\begin{aligned}
&\min
\sum_{k=1}^\wp \int_{t_\circ}^{t_\circ^N} \Big\| y^k(t)-\sum_{i=1}^\mathfrak l {\langle y^k(t),\psi_i\rangle}_X\,\psi_i\Big\|_X^2\,\mathrm dt\\
&\hspace{0.5mm}\text{s.t. } \{\psi_i\}_{i=1}^\mathfrak l\subset X\text{ and }{\langle\psi_i,\psi_j\rangle}_X=\delta_{ij},~1 \le
i,j \le \mathfrak l,
\end{aligned}
\right.
\end{equation}
where the symbol $\delta_{ij}$ denotes the Kronecker symbol satisfying $\delta_{ii}=1$ and $\delta_{ij}=0$ for $i\neq j$. An optimal solution $\{\bar\psi_i\}_{i=1}^\mathfrak l$ to \eqref{Eq:PODMin} is called a {\em POD basis of rank $\mathfrak l$}. The solution to \eqref{Eq:PODMin} is given by the next theorem. For its proof we refer the reader to \cite[Theorem~2.13]{GV13}.

\begin{theorem}
\label{Th:PODMin}
Let $X$ be a separable real Hilbert space and $y_1^k,\ldots,y_n^k\in X$ be given snapshots for $1\le k\le\wp$. Define the linear operator $\mathcal R:X\to X$ as follows:
\begin{equation}
\label{GVLuminy-Eq:OperatorR}
\mathcal R\psi=\sum_{k=1}^\wp\int_{t_\circ}^{t_\circ^N} {\langle \psi,y^k(t)\rangle}_X\,y^k(t)\,\mathrm dt\quad\text{for }\psi \in X.
\end{equation}
Then, $\mathcal R$ is a compact, nonnegative and symmetric operator. Suppose that $\{\bar\lambda_i\}_{i\in\mathbb N}$ and $\{\bar\psi_i\}_{i\in\mathbb N}$ denote the nonnegative eigenvalues and associated orthonormal eigenfunctions of $\mathcal R$ satisfying
\begin{equation}
\label{Eq:OperatorREigV}
\mathcal R\bar\psi_i=\bar\lambda_i\bar\psi_i,\quad\bar\lambda_1\ge\ldots\ge\bar\lambda_d>\bar\lambda_{d+1}=\ldots =0,\quad\bar\lambda_i \to 0 \text{ as } i \to \infty.
\end{equation}
Then, for every $\mathfrak l\le d$ the first $\mathfrak l$ eigenfunctions $\{\bar\psi_i\}_{i=1}^\mathfrak l$ solve \eqref{Eq:PODMin}. Moreover, the value of the cost evaluated at the optimal solution $\{\bar\psi_i\}_{i=1}^\mathfrak l$ satisfies
\begin{equation}
\label{Eq:PODError}
\mathcal E(\mathfrak l)=\sum_{k=1}^\wp\int_{t_\circ}^{t_\circ^N}\Big\|y^k(t)-\sum_{i=1}^\mathfrak l{\langle y^k(t),\bar\psi_i\rangle}_X\,\bar\psi_i\Big\|_X^2\,\mathrm dt=\sum_{i=\mathfrak l+1}^d\bar\lambda_i.
\end{equation}
\end{theorem}

\begin{remark}
\rm
In real computations, we do not have the whole trajectories $y^k(t)$ at hand f.a.a. $t\in [t_\circ,t_\circ^N]$ and for $1\le k\le\wp$. Moreover, the space $X$ has to be discretized as well. In this case, a discrete version of the POD method should be utilized; see, e.g., \cite{GV13}.\hfill$\Diamond$
\end{remark}

\subsubsection{The Galerkin POD scheme for the state equation}
\label{Section:4.2.2}

Suppose that $(t_\circ,y_\circ)\in\mathbb R^+_0\times H$ and $t_\circ^N=t_\circ+N\Delta t$ with prediction horizon $N\Delta t>0$. For given fixed control $u\in\mathbb U_{ad}^N(t_\circ)$ we consider the state equation $e(y,u)=0\in \mathbb Z^N(t_\circ)'$, i.e., $y$ satisfies \eqref{StateEquation}. Let us turn to a POD discretization of \eqref{StateEquation}. To keep the notation simple we apply only a spatial discretization with POD basis functions, but no time integration by, e.g., the implicit Euler method. In this section we distinguish two choices for $X$: $X=H$ and $X=V$. We choose the snapshots $y^1=y$ and $y^2=y_t$, i.e., we set $\wp=2$. By Proposition~\ref{Pro:StateEq} the snapshots $y^k$, $k=1,\ldots,\wp$, belong to $L^2(t_\circ,t_\circ^N;V)$. According to \eqref{Eq:OperatorREigV} let us introduce the following notations:
\begin{align*}
\mathcal R_V\psi&=\sum_{k=1}^\wp\int_{t_\circ}^{t_\circ^N} {\langle\psi,y^k(t) \rangle}_V\,y^k(t)\,\mathrm dt&&\text{for } \psi \in V,\\
\mathcal R_H\psi&=\sum_{k=1}^\wp\int_{t_\circ}^{t_\circ^N} {\langle\psi,y^k(t) \rangle}_H\,y^k(t)\,\mathrm dt&&\text{for } \psi \in H.
\end{align*}
To distinguish the two choices for the Hilbert space $X$ we denote by the sequence $\{(\lambda_i^V,\psi_i^V)\}_{i\in\mathbb N}\subset \mathbb R^+_0\times V$ the eigenvalue decomposition for $X=V$, i.e., we have
\[
\mathcal R_V\psi_i^V=\lambda_i^V\psi_i^V\quad \text{for all }i\in\mathbb N.
\]
Furthermore, let $\{(\lambda_i^H,\psi_i^H)\}_{i\in\mathbb N}\subset \mathbb R^+_0\times H$ in satisfy
\[
\mathcal R_H\psi_i^H=\lambda_i^H\psi_i^H\quad \text{for all }i\in\mathbb N.
\]
Then, $d=\dim \mathcal R_V(V)=\dim \mathcal R_H(H)\le\infty$; see \cite{Sin13}. The next result -- also taken from \cite{Sin13} -- ensures that the POD basis $\{\psi_i^H\}_{i=1}^\mathfrak l$ of rank $\mathfrak l$ build a subset of the test space $V$.

\begin{lemma}
\label{Lemma3.2.1}
Suppose that the snapshots $\{y^k\}_{k=1}^\wp$ belong to $L^2(t_\circ,t_\circ^N;V)$. Then, we have $\psi^H_i\in V$ for $i=1,\ldots,d$.
\end{lemma}

Let us define the two POD subspaces
\[
V^\mathfrak l=\mathrm{span}\,\big\{\psi_1^V,\ldots,\psi_\mathfrak l^V\big\}\subset V,\quad H^\mathfrak l=\mathrm{span}\,\big\{\psi_1^H,\ldots,\psi_\mathfrak l^H\big\}\subset V\subset H,
\]
where $H^\mathfrak l\subset V$ follows from Lemma~\ref{Lemma3.2.1}. Moreover, we introduce the orthogonal projection operators $\mathcal P^\mathfrak l_H:V\to H^\mathfrak l\subset V$ and $\mathcal P^\mathfrak l_V:V\to V^\mathfrak l\subset V$ as follows:
\begin{equation}
\label{GVLuminy:Eq3.2.12}
\begin{aligned}
&v^\mathfrak l=\mathcal P^\mathfrak l_H\varphi\text{ for any }\varphi\in V&&\text{iff } v^\mathfrak l\text{ solves } \min_{w^\mathfrak l\in H^\mathfrak l}{\|\varphi-w^\mathfrak l\|}_V,\\
&v^\mathfrak l=\mathcal P^\mathfrak l_V\varphi\text{ for any }\varphi\in V&&\text{iff } v^\mathfrak l\text{ solves } \min_{w^\mathfrak l\in V^\mathfrak l}{\|\varphi-w^\mathfrak l\|}_V.
\end{aligned}
\end{equation}
It follows from the first-order optimality conditions for \eqref{GVLuminy:Eq3.2.12} that $v^\mathfrak l=\mathcal P^\mathfrak l_H\varphi$ satisfies
\begin{equation}
\label{GVLuminy:Eq3.2.13}
{\langle v^\mathfrak l,\psi_i^H\rangle}_V={\langle\varphi,\psi_i^H\rangle}_V,\quad 1\le i\le \mathfrak l.
\end{equation}
Writing $v^\mathfrak l\in H^\mathfrak l$ in the form $v^\mathfrak l=\sum_{j=1}^\mathfrak l \mathrm v_j^\mathfrak l\psi_j^H$ we derive from \eqref{GVLuminy:Eq3.2.13} that the vector $\mathrm v^\mathfrak l=(\mathrm v_1^\mathfrak l,\ldots,\mathrm v_\mathfrak l^\mathfrak l)^\top\in\mathbb R^\mathfrak l$ satisfies the linear system
\begin{equation}
\label{Eq3.2.14}
\sum_{j=1}^\mathfrak l{\langle\psi_j^H,\psi_i^H\rangle}_V\,\mathrm v_j^\mathfrak l={\langle\varphi,\psi_i^H\rangle}_V,\quad 1\le i\le \mathfrak l.
\end{equation}
Summarizing, $v^\mathfrak l=\mathcal P^\mathfrak l_H\varphi\in H^\mathfrak l$ is given by the expansion $\sum_{j=1}^\mathfrak l \mathrm v_j^\mathfrak l\psi_j^H$, where the coefficients $\{\mathrm v_j^\mathfrak l\}_{j=1}^\mathfrak l$ satisfy the linear system \eqref{Eq3.2.14}. For the operator $\mathcal P^\mathfrak l_V:V\to V^\mathfrak l$ we have the explicit representation
\begin{equation}
\label{Proj:POD}
\mathcal P^\mathfrak l_V\varphi=\sum_{i=1}^\mathfrak l \langle\varphi,\psi_i\rangle_V\,\psi_i\quad\text{for }\varphi\in V. 
\end{equation}
We conclude from \eqref{Eq:PODError} that
\begin{equation}
\label{Eq:RhoError-1}
\sum_{k=1}^\wp \int_{t_\circ}^{t_\circ^N} {\|y^k(t)-\mathcal P_V^\mathfrak l y^k(t)\|}_V^2\,\mathrm dt=\sum_{i=\mathfrak l+1}^d\lambda_i^V.
\end{equation}

Let us define the linear space $X^\mathfrak l\subset V$ as
\[
X^\mathfrak l=\mathrm{span}\,\big\{\psi_1,\ldots,\psi_\mathfrak l\big\},
\]
where $\psi_i=\psi_i^V$ in case of $X=V$ and $\psi_i=\psi_i^H$ in case of $X=H$. Hence, $X^\mathfrak l=V^\mathfrak l$ and $X^\mathfrak l=H^\mathfrak l$ for $X=V$ and $X=H$, respectively. Now, a POD Galerkin scheme for \eqref{StateEquation} is given as follows: find $y^\mathfrak l(t)\in X^\mathfrak l$ f.a.a. $t\in[t_\circ,t_\circ^N]$ satisfying
\begin{equation}
\label{Eq:PODGalState}
\begin{aligned}
&\frac{\mathrm d}{\mathrm dt}\,{\langle y^\mathfrak l(t),\psi\rangle}_H+\int_\Omega\theta y^\mathfrak l_x(t)\psi'+\big(y^\mathfrak l_x(t)+\rho(y^\mathfrak l(t)^3-y^\mathfrak l(t))\big)\psi\,\mathrm dx\\
&\hspace{50mm}=\int_\Omega u(t)\psi\,\mathrm dx\quad\text{f.a.a. }t\in(t_\circ,t_\circ^N],\\
&{\langle y^\mathfrak l(t_\circ),\psi\rangle}_H={\langle y_\circ,\psi\rangle}_H
\end{aligned}
\end{equation}
for all $\psi\in X^\mathfrak l$. It follows by similar arguments as in the proof of Proposition~\ref{Pro:StateEq} that there exists a unique solution to \eqref{Eq:PODGalState}. If $y_\circ\in L^\infty(Q^N)$ holds, $y^\mathfrak l$ satisfies the a-priori estimate
\begin{equation}
\label{Eq:LinfEstPOD}
{\|y^\mathfrak l\|}_{\mathbb Y^N(t_\circ)}+{\|y^\mathfrak l\|}_{L^\infty(Q^N)}\le C\big({\|y_\circ\|}_{L^\infty(\Omega)}+{\|u\|}_{\mathbb U^N(t_\circ)}\big),
\end{equation}
where the constant $C>0$ is independent of $\mathfrak l$ and $y_\circ$. Let $\mathcal P^\mathfrak l$ denote $\mathcal P_V^\mathfrak l$ in case of $X=V$ and $\mathcal P_H^\mathfrak l$ in case of $X=H$. The next result is proved in \ref{AppA}.

\begin{theorem}
\label{Pro:Apriori}
Suppose that $(t_\circ,y_\circ)\in\mathbb R^+_0\times L^\infty(\Omega)$, $t_\circ^N=t_\circ+N\Delta$ with prediction horizon $N\Delta t>0$. Further, let $u\in \mathbb U_{ad}^N(t_\circ)$ be a fixed control input. By $y$ and $y^\mathfrak l$ we denote the unique solution to \eqref{StateEquation} and \eqref{Eq:PODGalState}, respectively, where the POD basis of rank $\mathfrak l$ is computed by choosing $\wp=2$, $y^1=y$ and $y^2=y_t$. Then,
\[
{\|y-y^\mathfrak l\|}_{\mathbb Y^N(t_\circ)}^2\le C\cdot\left\{
\begin{aligned}
&{\|y^\mathfrak l(t_\circ)-\mathcal P^\mathfrak l_V y_\circ\|}_H^2+\sum_{i=\mathfrak l+1}^d\lambda_i^V,&&X=V,\\
&{\|y^\mathfrak l(t_\circ)-\mathcal P^\mathfrak l_H y_\circ\|}_H^2+\sum_{i=\mathfrak l+1}^d\lambda_i^H{\|\psi_i^H-\mathcal P_H^\mathfrak l\psi_i^H\|}_V^2,&&X=H
\end{aligned}
\right.
\]
for a $C>0$ which is independent of $\mathfrak l$. In particular, $\lim_{\mathfrak l\to\infty}\|y-y^\mathfrak l\|_{\mathbb Y^N(t_\circ)}=0$.
\end{theorem}

\subsubsection{The Galerkin POD scheme for the optimality system}
\label{Section:4.2.3}

Suppose that we have computed a POD basis $\{\psi_i\}_{i=1}^{\mathfrak l}$ of rank $\mathfrak l$ by choosing $X=H$ or $X=V$. Suppose that for $u\in\mathbb U_{ad}^N(t_\circ)$ the function $y^\mathfrak l$ is the POD Galerkin solution to \eqref{Eq:PODGalState}. Then the POD Galerkin scheme for the adjoint equation \eqref{DualEq} is given as follows: find $p^\mathfrak l\in X^\mathfrak l=\mathrm{span}\,\{\psi_1,\ldots,\psi_\mathfrak l\}$ f.a.a. $t\in[t_\circ,t_\circ^N]$ satisfying
\begin{equation}
\label{PODDualEq}
\begin{aligned}
&-\frac{\mathrm d}{\mathrm dt}\,{\langle p^\mathfrak l(t),\psi\rangle}_H+\int_\Omega\theta p^\mathfrak l_x(t)\psi'-\big(p^\mathfrak l_x(t)+\rho(1-3y^\mathfrak l(t)^2)\big)p^\mathfrak l(t)\psi\,\mathrm dx\\
&\hspace{36mm}=\int_\Omega \big(y_d-y^\mathfrak l(t)\big)\psi\,\mathrm dx=0\quad\text{f.a.a. }t\in[t_\circ,t_\circ^N),\\
&{\langle p^\mathfrak l(t_\circ^N),\psi\rangle}_H=0
\end{aligned}
\end{equation}
for all $\psi\in X^\mathfrak l$. A-priori error estimates for the POD solution $p^\mathfrak l$ to \eqref{PODDualEq} can be derived by variational arguments; compare \cite{SS13} and \cite[Theorem~4.15]{GV13}. If $p^\mathfrak l$ is computed, we can derive a POD approximation for the variational inequality \eqref{VarIneq}:
\begin{equation}
\label{PODVarIneq}
\int_{t_\circ}^{t_\circ^N}\int_\Omega(\lambda u-p^\mathfrak l)(\tilde u-u)\,\mathrm dx\mathrm dt\ge 0\quad \text{for all }\tilde u\in\mathbb U^N_{ad}(t_\circ).
\end{equation}
Summarizing, a POD suboptimal solution $\bar x^{N,\mathfrak l}=(\bar y^{N,\mathfrak l},\bar u^{N,\mathfrak l})\in\mathbb X_{ad}^N(t_\circ)$ to \eqref{P-DynSys-FinHor} satisfies together with the associated Lagrange multiplier $\bar p^{N,\mathfrak l}\in\mathbb Y_1^N(t_\circ)$ the coupled system \eqref{Eq:PODGalState}, \eqref{PODDualEq} and \eqref{PODVarIneq}. The POD approximation of the finite horizon quadratic cost functional \eqref{PN-OptPro} reads
\[
\hat J^{N,\mathfrak l}(u;t_\circ,y_\circ)=\int_{t_\circ}^{t_\circ^N} \ell\big(y^\mathfrak l_{[u,t_\circ,y_\circ]}(t),u(t)\big)\,\mathrm dt,
\]
where $y^\mathfrak l_{[u,t_\circ,y_\circ]}$ is the solution to \eqref{Eq:PODGalState}. In Algorithm~\ref{Alg:POD-NMPC} we set up the POD discretization for Algorithm~\ref{Alg:NMPC}.
\begin{algorithm}[htbp]
\caption{(POD-NMPC algorithm)}
\label{Alg:POD-NMPC}
\begin{algorithmic}[1]
\REQUIRE time step $\Delta t>0$, finite control horizon $N\in\mathbb N$, weighting parameter $\lambda>0$, POD tolerance $\tau_{pod}>0$.
\STATE Compute a POD basis $\{\psi_i\}_{i=1}^\mathfrak l$ satisfying \eqref{Eq:PODError} with $\mathcal E(\mathfrak l)\le \tau_{pod}$.
\FOR{$n=0,1,2,\ldots$}
\STATE Measure the state $y(t_n)\in V$ of the system at $t_n=n\Delta t$.
\STATE Set $t_\circ=t_n=n\Delta t$, $y_\circ=y(t_n)$ and compute a global solution to 
\begin{equation}
\tag{$\mathbf{P^{\boldsymbol{N,\mathfrak l}}(}t_\circ\mathbf{)}$}
\label{P-DynSysPOD-FinHor}
\min\hat J^{N,\mathfrak l}(u^\mathfrak l;t_\circ,y^\mathfrak l_\circ)\quad\text{s.t.}\quad u^\mathfrak l\in\mathbb U_{ad}^N(t_\circ).
\end{equation}
We denote the optimal control by $\bar u^{N,\mathfrak l}$ and the optimal state by $\bar y^{N,\mathfrak l}$.
\STATE Define the NMPC feedback value $\mu^{N,\mathfrak l}(t;t_\circ,y_\circ)=\bar u^{N,\mathfrak l}(t)$ and use this control to compute the associated state $y=y_{[\mu^{N,\mathfrak l}(\cdot),t_\circ,y_\circ]}$ by solving \eqref{DynamicalSystem} on $[t_\circ,t_\circ+\Delta t]$.
\ENDFOR
\end{algorithmic}
\end{algorithm}
Due to our POD reduced-order approach an optimal solution to \eqref{P-DynSysPOD-FinHor} can be computed much faster than the one to \eqref{P-DynSys-FinHor}. In the next subsection we address the question, how the suboptimality of the control influences the asymptotic stability.

\subsection{Asymptotic stability for the POD-MPC algorithm}
\label{Section:4.3}

In this subsection we present the main results of this paper. We give sufficient conditions that Algorithm~\ref{Alg:POD-NMPC} gives a stabilizing feedback control for the reduced-order model. Due to Definition~\ref{De:ExpSta} we have to find an admissible control $u\in\mathbb U^N(t_\circ)$ for any $N\in\mathbb N$ so that the solution to \eqref{DynamicalSystem} satisfies \eqref{exp:con}.

In \eqref{RunCost} we have introduced our running quadratic cost. As in Section~\ref{Section3.2} we choose $y_d=y_*=0$. Suppose that $y^\mathfrak l$ is the reduced-order solution to \eqref{Eq:PODGalState} for the control $u^\mathfrak l=-Ky^\mathfrak l$. If $K$ satisfies appropriate bounds (see Remark~\ref{Rem:Kest}-2)), we can ensure that $u^\mathfrak l\in\mathbb U^N_{ad}(t_\circ)$ holds. Analogously to \eqref{AprioriEst} ands \eqref{stab:prp} we find
\begin{equation}
\label{stab:PODprp}
{\|y^\mathfrak l(t)\|}_H^2\le \sigma(K)^{t-t_\circ}\,{\|y_\circ\|}_H^2\quad\text{f.a.a. }t\ge t_\circ
\end{equation}
and
\begin{equation}
\label{Eq:Est1}
\ell\big(y^\mathfrak l(t),u^\mathfrak l(t)\big)\le\frac{C(K)}{2}\,{\|y^\mathfrak l(t)\|}_H^2.
\end{equation}
with the same constants $C(K)$ and $\sigma(K)$ as in \eqref{C(K)}. Let $y_{[u^\mathfrak l,t_\circ,y_\circ]}$ be the (full-order) solution to \eqref{Eq:PODGalState} for the same admissible control law $u=u^\mathfrak l$. Utilizing the Cauchy-Schwarz inequality we get
\begin{equation}
\label{Eq:Est2}
\begin{aligned}
\ell\big(y_{[u^\mathfrak l,t_\circ,y_\circ]}(t),u^\mathfrak l(t)\big)
&\le\frac{1}{2}\,{\|y_{[u^\mathfrak l,t_\circ,y_\circ]}(t)-y^\mathfrak l(t)\|}_H^2+\ell\big(y^\mathfrak l(t),u^\mathfrak l(t)\big)\\
&\qquad+{\|y_{[u^\mathfrak l,t_\circ,y_\circ]}(t)-y^\mathfrak l(t)\|}_H{\|y^\mathfrak l(t)\|}_H.
\end{aligned}
\end{equation}
If $y_\circ\neq0$ holds, we infer that $\|y^\mathfrak l(t)\|_H$ is positive for all $t\in[t_\circ,t_\circ^N]$. Then, we conclude from \eqref{Eq:Est1}, \eqref{Eq:Est2} and \eqref{stab:PODprp} that the exponential controllability condition \eqref{exp:con} holds for the admissible control law $u^\mathfrak l=-Ky^\mathfrak l$:
\begin{align*}
\ell\big(y_{[u^\mathfrak l,t_\circ,y_\circ]}(t),u^\mathfrak l(t)\big)&\le\frac{1}{2}\,\bigg(\mathsf{Err}(t;\mathfrak l)^2+C(K)+2\mathsf{Err}(t;\mathfrak l)\bigg)\,{\|y^{\mathfrak l}(t)\|}_H^2\\
&\le \frac{1}{2}\,C^\mathfrak l(K)\,\sigma(K)^{t-t_\circ}\,{\|y_\circ\|}_H^2=C^\mathfrak l(K)\,\sigma(K)^{t-t_\circ}\,\ell_*(y_\circ)
\end{align*}
with the error term
\begin{equation}
\label{Eq:ErrTerm}
\mathsf{Err}(t;\mathfrak l)=\frac{{\|y_{[u^\mathfrak l,t_\circ,y_\circ]}(t)-y^\mathfrak l(t)\|}_H}{{\|y^\mathfrak l(t)\|}_H}
\end{equation}
and the constant
\begin{equation}
\label{Cl}
C^\mathfrak l(K)=C(K)+2\mathsf{Err}(t;\mathfrak l)+\mathsf{Err}(t;\mathfrak l)^2\ge C(K).
\end{equation}
Thus, the constant $C^\mathfrak l(K)$ takes into account the approximation made by the POD reduced-order model. In the following theorem we provide an explicit formula for the scalar $\alpha^{N,\mathfrak l}$ which appears in the relaxed DPP. The notation $\alpha^{N,\mathfrak l}$ intends to stress that we are working with POD surrogate model. We summarize our result in the following theorem.

\begin{theorem}\label{Th:alphal}
Let the constant $C^\mathfrak l$ be given by \eqref{Cl} and $N\Delta t$ denote the finite prediction horizon with $N\in\mathbb N$ and $\Delta t>0$. Then the parameter $\alpha^{N,\mathfrak l}$ is given by the explicit formula:
\begin{equation}
\label{alphal}
\alpha^{N,\mathfrak l}(K)=1-\frac{\big(\eta^\mathfrak l_N(K)-1\big)\prod_{i=2}^N\big(\eta^\mathfrak l_i(K)-1\big)}{\prod_{i=2}^N\eta^\mathfrak l_i(K)-\prod_{i=2}^N\big(\eta^\mathfrak l_i(K)-1\big)}
\end{equation}
with $\eta^\mathfrak l_i(K)=C^\mathfrak l(K)(1-\sigma^i(K))/(1-\sigma(K))$  and $\sigma(K)$ as in \eqref{C(K)}.

\end{theorem}

\begin{remark}
\rm
\begin{enumerate}
\item [1)] If $\mathsf{Err}(t;\mathfrak l)$ is small, Theorem \ref{Th:alphal} informs we can compute the constant $\alpha^{N,\mathfrak l}\approx\alpha^N$ basically in the same way of the full-model, replacing the constants $C$, $\eta$ with $C^\mathfrak l$, $\eta^\mathfrak l$, respectively, taking into account the POD reduced-order modelling. Then, \eqref{rel:est} implies that a suboptimality estimate holds approximately; see Remark~\ref{rmk_imp}. To obtain the minimal horizon which ensures the asymptotic stability of the POD-NMPC scheme we maximize \eqref{alphal} according to the constraints $\alpha^{N,\mathfrak l}>0$, $K>\max(0,\rho-\theta/C_V)$ and to the constraints in Table \ref{K:cond1}.
\item [2)] Due to \eqref{stab:PODprp} and $u^\mathfrak l=-Ky^\mathfrak l$  the norm $\|u^\mathfrak l(t)\|_H$ is bounded independent of $\mathfrak l$. By Theorem~\ref{Pro:Apriori} and \eqref{Embed-2} we have $\lim_{\mathfrak l\to\infty}\|y_{[u^\mathfrak l,t_\circ,y_\circ]}(t)-y^\mathfrak l(t)\|_H=0$ holds for all $t\in[t_\circ,t_\circ^N]$. Thus, if we choose $\mathfrak l$ sufficiently large we can ensure that $\mathsf{Err}(t;\mathfrak l)$ is small enough provided the denominator satisfies $\|y^\mathfrak l(t)\|_H\ge C_*$ with a positive constant $C_*$ which is independent of $\mathfrak l$.
\item [3)] In Algorithm~\ref{Alg:POD-NMPC} we compute the control law $\bar u^{N,\mathfrak l}$ instead of $-Ky^\mathfrak l$. Therefore, one can replace $\mathsf{Err}(t;\mathfrak l)$ by
\[
\widetilde{\mathsf{Err}}(t;\mathfrak l)=\frac{{\|y_{[\bar u^{N,\mathfrak l}(\cdot),t_\circ,y_\circ]}(t)-y^{N,\mathfrak l}(t)\|}_H}{{\|y^{N,\mathfrak l}(t)\|}_H}
\]
that can be evaluated easily, since $y^{N,\mathfrak l}(t)$ and $y_{[\bar u^{N,\mathfrak l},t_\circ,y_\circ]}$ are known from Algorithm~\ref{Alg:POD-NMPC}, steps 4 and 5, respectively. It turns out that for our test examples both error terms lead to the same choices for the prediction horizon $N\in\mathbb N$, for the positive feedback factor $K$ and for the relaxation parameter $\alpha^{N,\mathfrak l}\in(0,1]$.\hfill$\Diamond$
\end{enumerate}
\end{remark}

\section{Numerical tests}
\label{Section5}
\setcounter{equation}{0}
\setcounter{theorem}{0}
\setcounter{table}{0}

This section presents numerical tests in order to show the performance of our proposed algorithm. All the numerical simulations reported in this paper have been made on a MacBook Pro with 1 CPU Intel Core i5 2.3 Ghz and 8GB RAM.

\subsection{The finite difference approximation for the state equation}

\begin{figure}[H]
\centering
\includegraphics[scale=0.38]{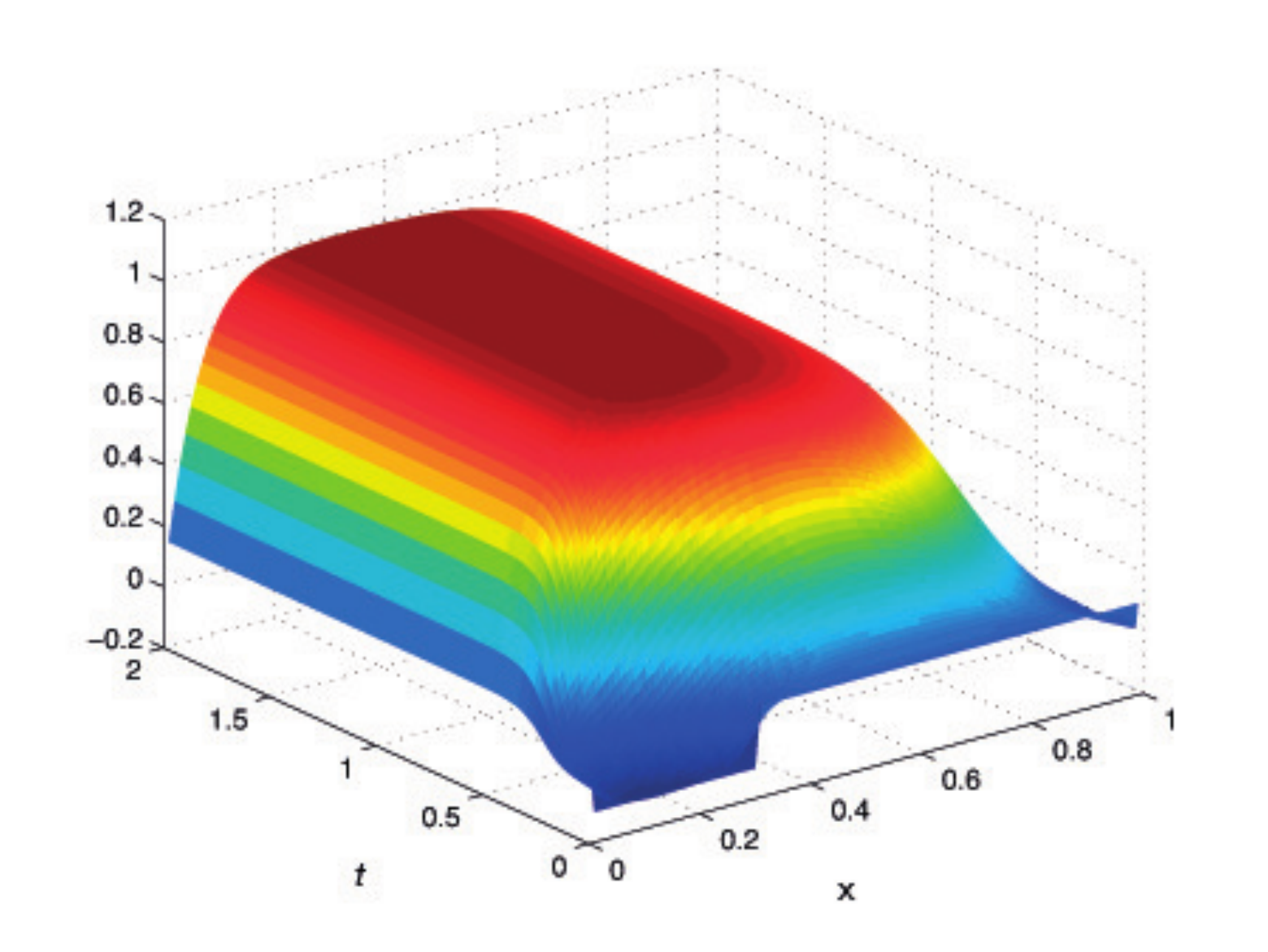}\hspace{-2mm}
\includegraphics[scale=0.38]{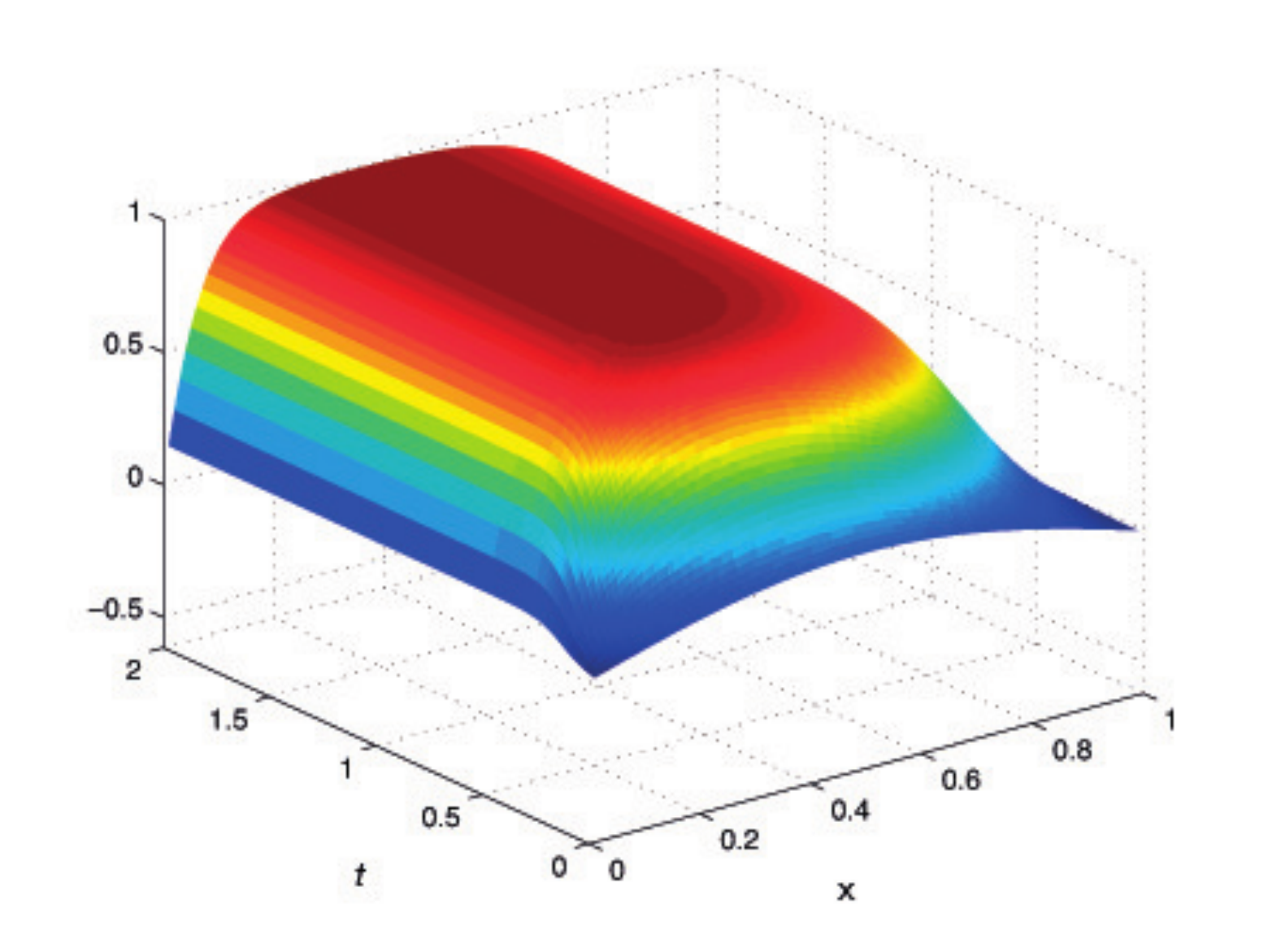}
\caption{FD state $y$ for $y_\circ=0.1\sgn(x-0.3)$ (left plot) and $y_\circ=0.2\sin{\pi x}$ (right plot) with $u=0$, $(\theta,\rho)=(0.1,11)$ and $\mathcal N=99$.}
\label{fig:1}
\end{figure}

For $\mathcal N\in\mathbb N$ we introduce an equidistant spatial grid in $\Omega$ by $x_i=i\Delta x$, $i=0,\ldots,\mathcal N+1$, with the step size $\Delta x=1/(\mathcal N+1)$. At $x_0=0$ and $x_{\mathcal N+1}=1$ the solution $y$ is known due to the boundary conditions \eqref{PDESystem}. Thus, we only compute approximations $y_i^h(t)$ for $y(t,x_i)$ with $1\le i\le\mathcal N$ and $t\in[t_\circ,t_f]$. We define the vector $y^h(t)=(y_1^h(t),\ldots,y_\mathcal N^h(t))^\top\in\mathbb R^\mathcal N$ of the unknowns. Analogously, we define $u^h=(u^h_1,\ldots,u^h_\mathcal N)^\top\in\mathbb R^\mathcal N$, where $u^h_i$ approximates $u(x_i,\cdot)$ for $1\le i\le\mathcal N$. Utilizing a classical second-order finite difference (FD) scheme and an implicit Euler method for the time integration we derive a discrete approximation of the parabolic problem. In Figure~\ref{fig:1} the discrete solutions are plotted for $\mathcal N=99$, for $t\in[0,2]$ and two different initial conditions.
As we see from Figure~\ref{fig:1}, the uncontrolled solutions do not tend to zero for $t\to\infty$, indeed it stabilizes at one.

\subsection{POD-NMPC experiments}\label{num_test}

\begin{figure}[htbp]
\centering
\includegraphics[scale=0.38]{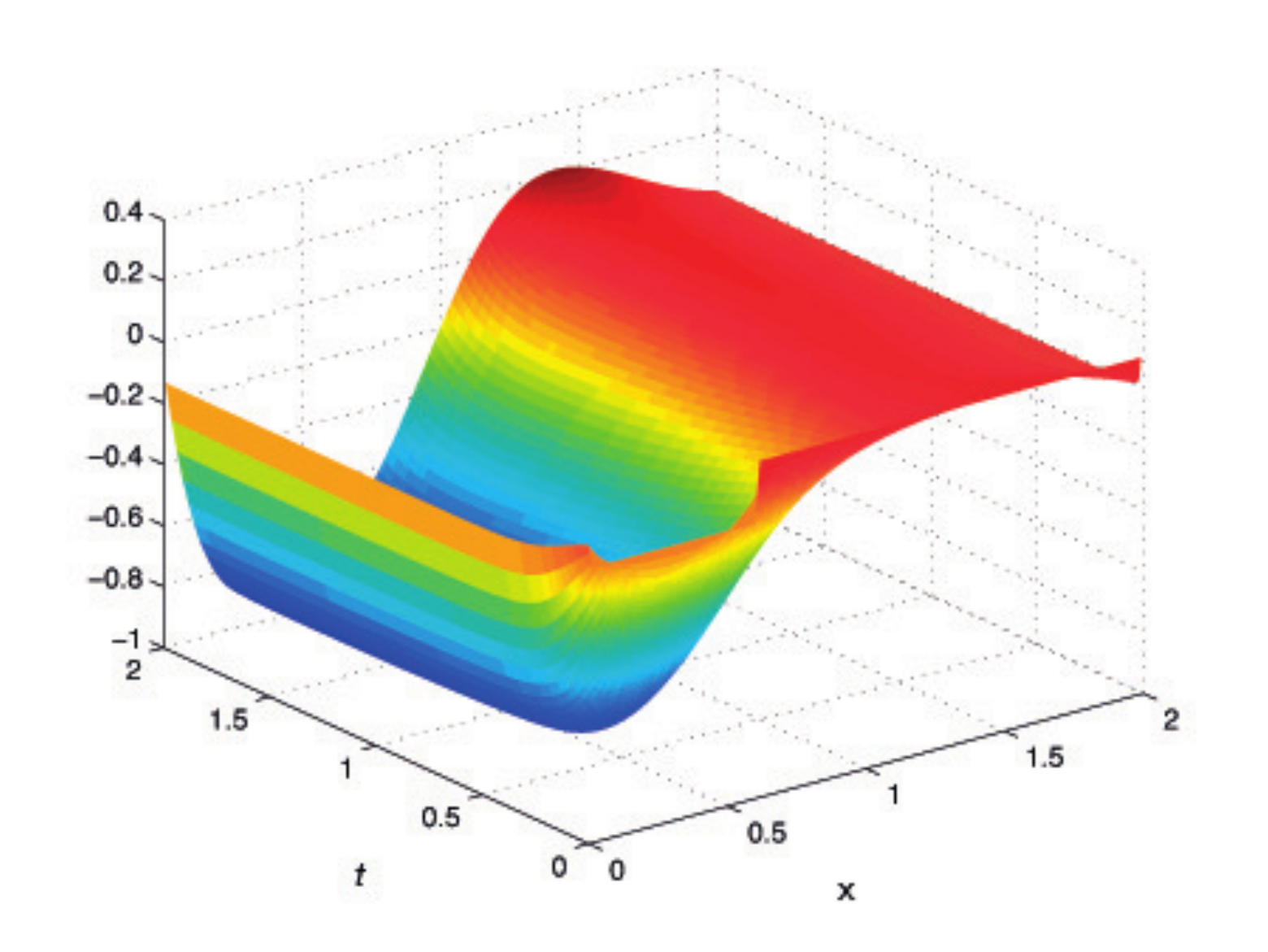}\hspace{-2mm}
\includegraphics[scale=0.38]{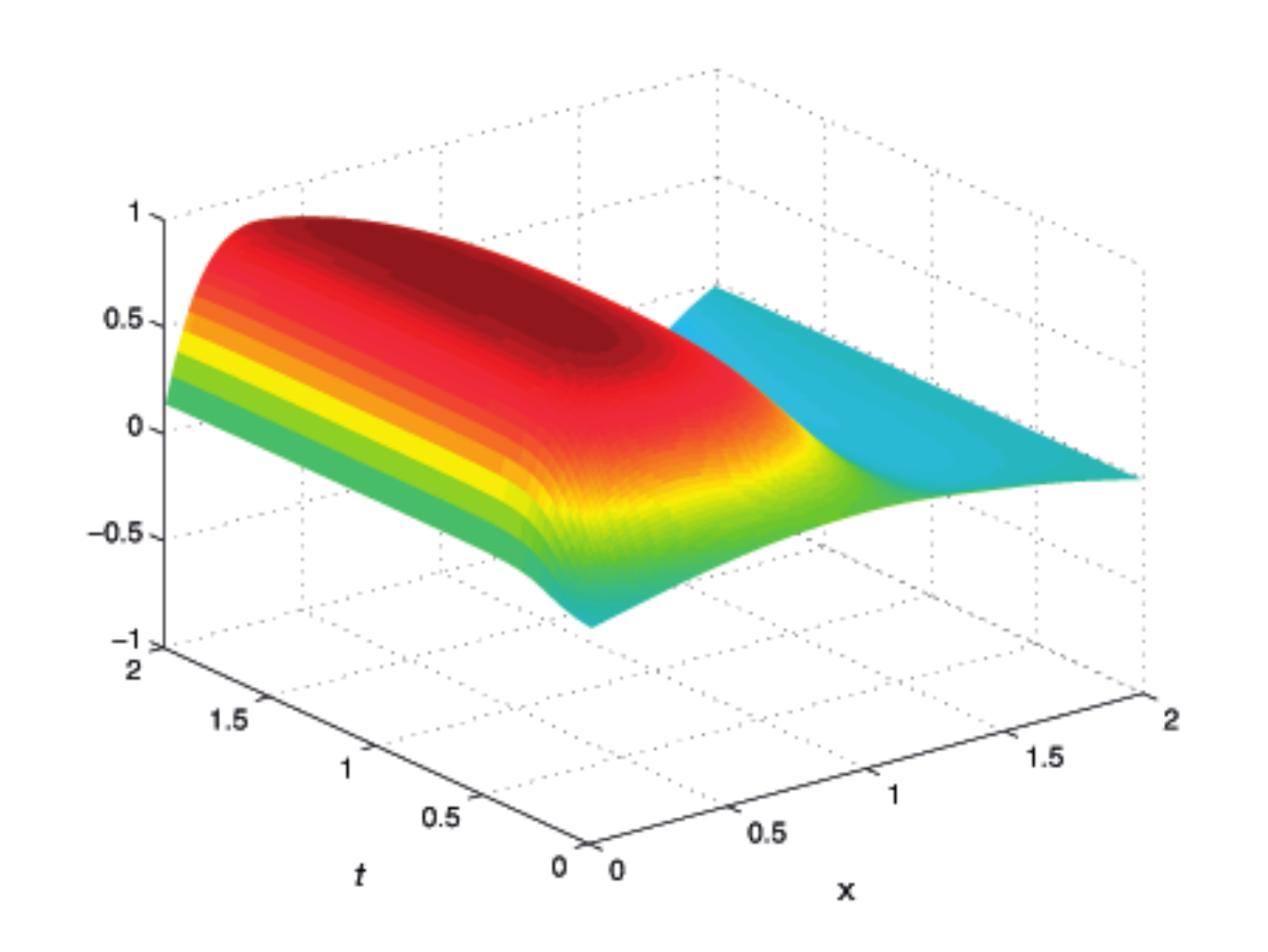}
\caption{Open-loop solution  $y$ for $y_\circ=0.1\sgn(x-0.3)$, $(\theta,\rho)=(0.1,11)$, $\mathcal N=99$, $t_f=2$ (left plot) and $y_\circ=0.2\sin{\pi x}$, $(\theta,\rho)=(0.1,11)$, $\mathcal N=99$, $t_f=2$ (right plot). }
\label{fig:2}
\end{figure}

In our numerical examples we choose $y_d\equiv0,$ i.e., we force the state to be close to zero, and $\lambda=0.01$ in \eqref{InfiniteHorizonCost}. A finite horizon open loop strategy does not steer the trajectory to the zero-equilibrium (see Figure~\ref{fig:2}). Therefore, stabilization is not guaranteed by the theory of asymptotic stability. Note that we are not dealing with terminal constraints and the terminal condition of the adjoint equation \eqref{DualEq} is zero.
In our tests, the snapshots are computed taking the uncontrolled system, e.g.  $u\equiv 0,$ in \eqref{PDESystem} and the correspondent adjoint equation \eqref{DualEq}. Several hints for the computation of the snapshots in the context of MPC are given in \cite{GU12}. The nonlinear term is reduced following the Discrete Empirical Interpolation Method (DEIM) which is a method that avoid the evaluation of the full model of the nonlinear part building new basis functions upon the nonlinear term; compare \cite{CS09} for more details. Note that, in our simulations, the optimal prediction horizon $N$ is always obtained from Theorem \ref{Th:alphal}.

\begin{run}[\em Unconstrained case with smooth initial data]
\label{Run1}
\rm
The parameters are presented in Table~\ref{par:1}.
\begin{table}[htbp]
\begin{center}
\begin{tabular}{cccccccccc}
\toprule
$T$ & $\Delta t$ & $\Delta x$ & $\theta$ & $\rho$ & $y_0(x)$ & $u_a$ & $u_b$ & $N$ & $K$\\
\midrule
0.5 & 0.01 & 0.01 & 1 & 11 & $0.2\sin(\pi x)$ & $-\infty$ & $\infty$ & 10 & 2.46\\
\bottomrule
\end{tabular}
\end{center}
\caption{Run~\ref{Run1}: Setting for the optimal control problem, minimal stabilizing horizon $N$ and feedback constant $K$.}
\label{par:1}
\end{table}
According to the computation of $\alpha^N$ in \eqref{alpha} related to the relaxed DPP, the minimal horizon that guarantes asymptotic stability is $N=10$. Even in the POD-NMPC scheme the asymptotic stability is achieved for $N=10$, provided that $\mathsf{Err}(t;\mathfrak l)\le10^{-3}$ for all $t\ge t_\circ$. Note that the horizon of the surrogate model is computed by \eqref{alphal}. In Figure~\ref{fig:test1} we show the controlled state trajectory computed by Algorithm~\ref{Alg:NMPC} taking $N=3$ and $N=10$.
\begin{figure}[htbp]
\includegraphics[height=43mm,width=40mm]{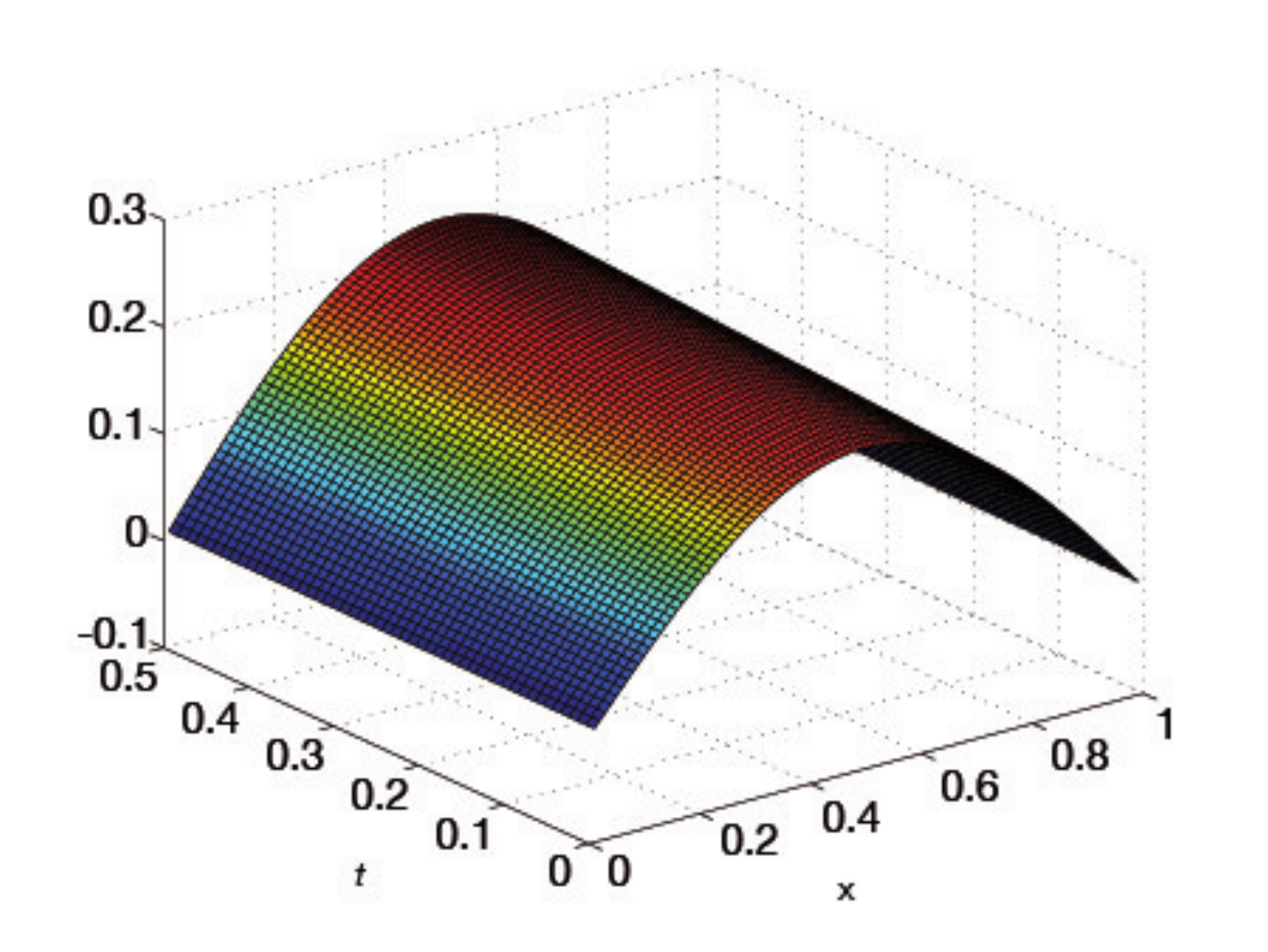}\hfill\includegraphics[height=43mm,width=40mm]{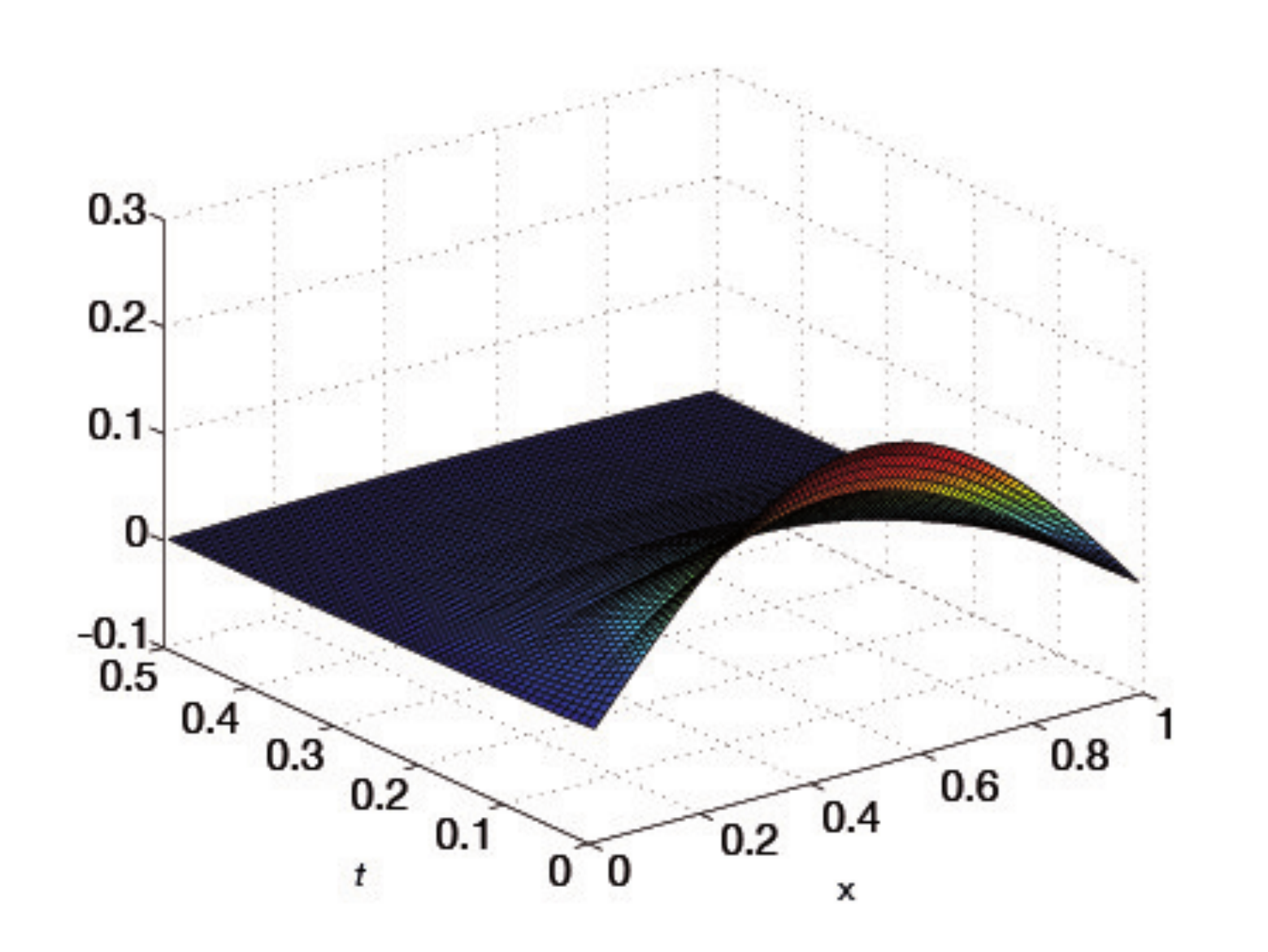}\hfill\includegraphics[height=43mm,width=40mm]{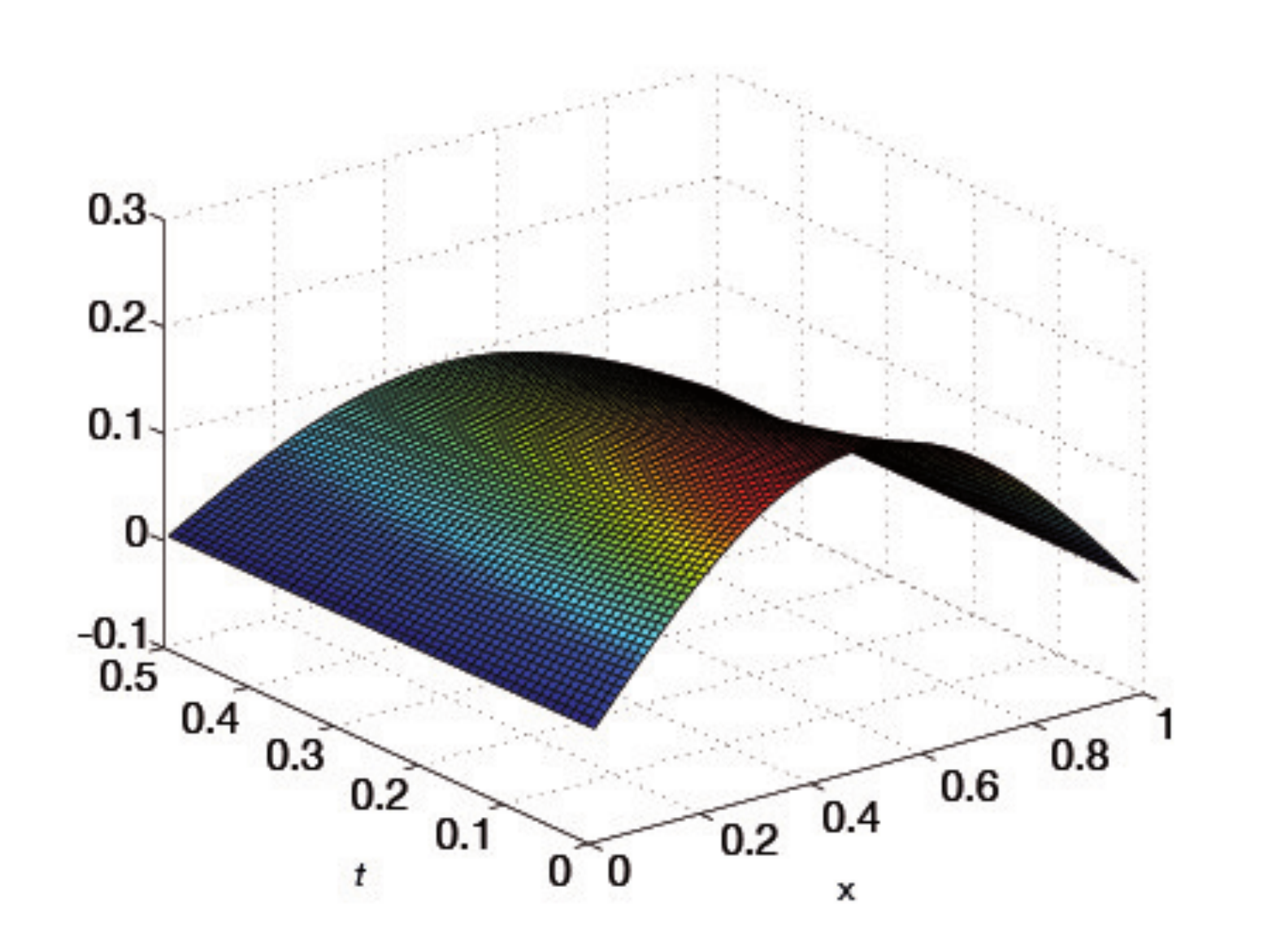}
\caption{Run~\ref{Run1}: NMPC state with $N=3$ (left plot), with $N=10$ (middle plot) and with $u=-Ky$ (right plot)}
\label{fig:test1}
\end{figure}
As we can see, we do not get a stabilizing feed\-back for $N=3$, whereas $N=10$ leads to a state trajectory which tends to zero for $t\to\infty$. Note that we plot the solution only on the time interval $[0, 0.5]$ in order to have a zoom of the solution. Further, in Figure~\ref{fig:test1} the solution related to $u=-Ky$ is presented. As we can see, the NMPC control stabilized to the origin very soon while the control law $u=-Ky$ requires a larger time horizon. This is due to the fact that MPC stabilizes in an optimal way, in contrast to the control law $u=-Ky$. In Table~\ref{tab:test1} we present the error in $L^2(t_\circ,T;H)$-norm considering the solution coming from the  Algorithm~\ref{Alg:NMPC} as the 'truth' solution (in our case the finite difference solution denoted by $y^{FD}$). The examples are computed with $\mathsf{Err}(t,\mathfrak l)\le 10^{-3}.$
\begin{table}[htbp]
\begin{center}
\begin{tabular}{lcccc}
\toprule
 & $\hat J$ & time & $K$ & $\|y^{FD}-y\|_{L^2(t_\circ,T;H)}$\\
\midrule
Solution with $u=-Ky$ & 0.0025 &  & 2.46 & 0.0145\\
Alg.~\ref{Alg:NMPC} & 0.0015 & 49s &   &\\
Alg.~\ref{Alg:POD-NMPC} ($\mathfrak l=13,\mathfrak l^{DEIM}=15$) & 0.0016 & \phantom{1}8s & & 0.0047\\
Alg.~\ref{Alg:POD-NMPC} ($\mathfrak l=3,\mathfrak l^{DEIM}=2$) & 0.0016 & \phantom{1}6s & & 0.0058\\
\bottomrule
\end{tabular}
\end{center}
\caption{Run~\ref{Run1}: Evaluation of the cost functional, CPU time, suboptimal solution.}
\label{tab:test1}
\end{table}
The CPU time for the full-model turns out to be 49 seconds, whereas the POD-suboptimal approximation with only three POD and two DEIM basis functions requires 6 seconds. We can easily observe an impressive speed up factor eight. Moreover the evaluation of the cost functional in the full model and the POD model provides very close values. We have not considered the CPU time in the suboptimal problem since it did not involve a real optimazion problem. As soon as we have computed $K$, within an offline stage, we directly approximate the equation with the control law $u=-Ky$.\hfill$\Diamond$
\end{run}

\begin{run}[\em Constrained case with smooth initial data]
\label{Run2}
\rm
In contrast to Run~\ref{Run1} we choose $u_a=-0.3$ and $u_b=0.$ As expected, the minimal horizon $N$ increases compared to Run~\ref{Run1}; see Table \ref{par:2}.
\begin{table}[htbp]
\begin{center}
\begin{tabular}{cccccccccc}
\toprule                                                       
$T$ & $\Delta t$ & $\Delta x$ & $\theta$ & $\rho$ & $y_0(x)$ & $u_a$ & $u_b$ & $N$  & $K$\\
\midrule
0.5 & 0.01 & 0.01 & 1 & 11 & $0.2\sin(\pi x)$ & $-0.3$ & $0$ & 14 & 1.50\\
\hline
\end{tabular}
\end{center}
\caption{Run~\ref{Run2}: Setting for the optimal control problem, minimal stabilizing horizon $N$ and feedback constant $K$.}
\label{par:2}
\end{table}
As one can see from Figure~\ref{pic:test2} the NMPC state with $N=14$ tends faster to zero than the state with $u=-Ky$.
\begin{figure}[htbp]
\includegraphics[height=43mm,width=40mm]{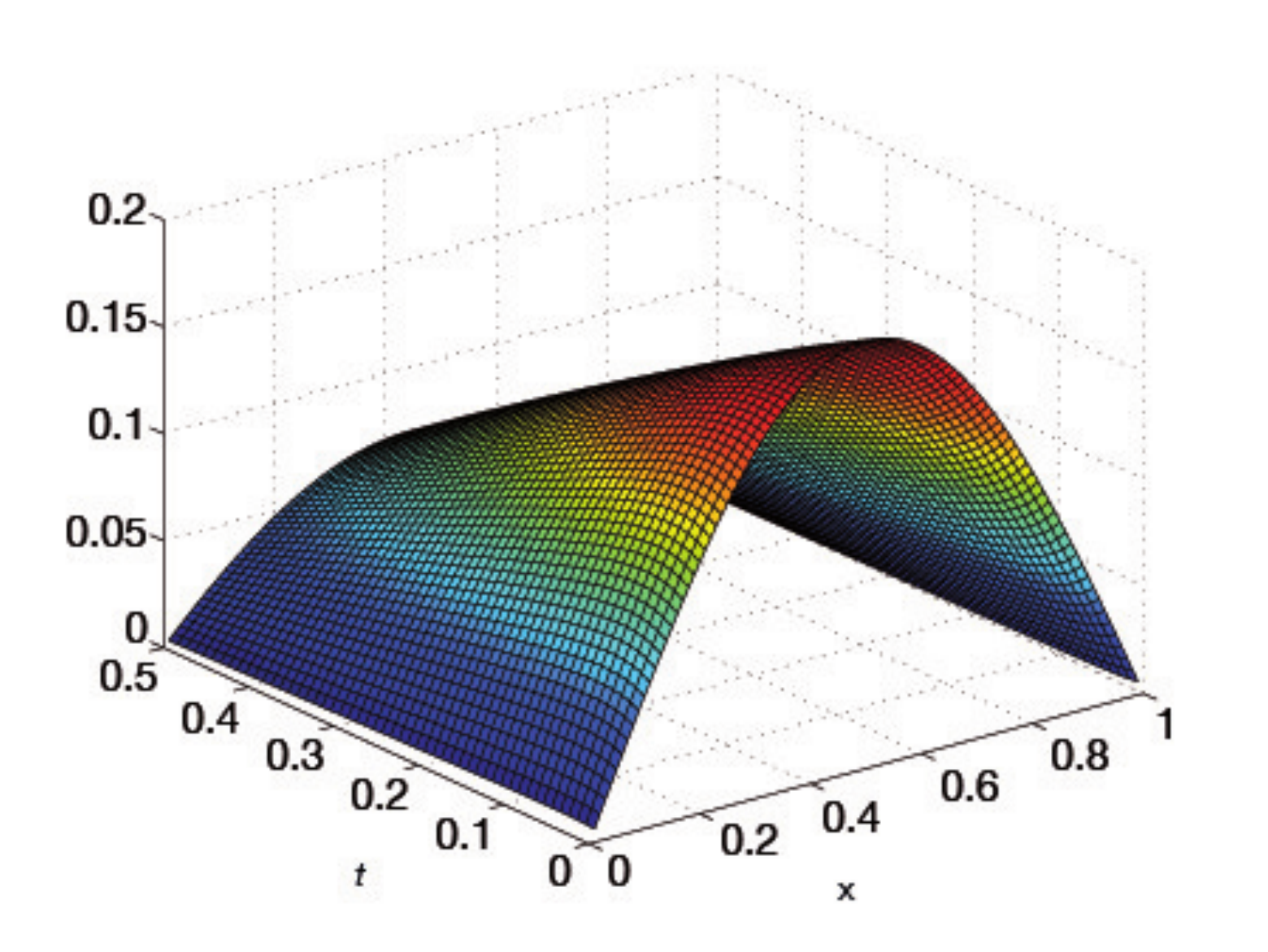}\hfill\includegraphics[height=43mm,width=40mm]{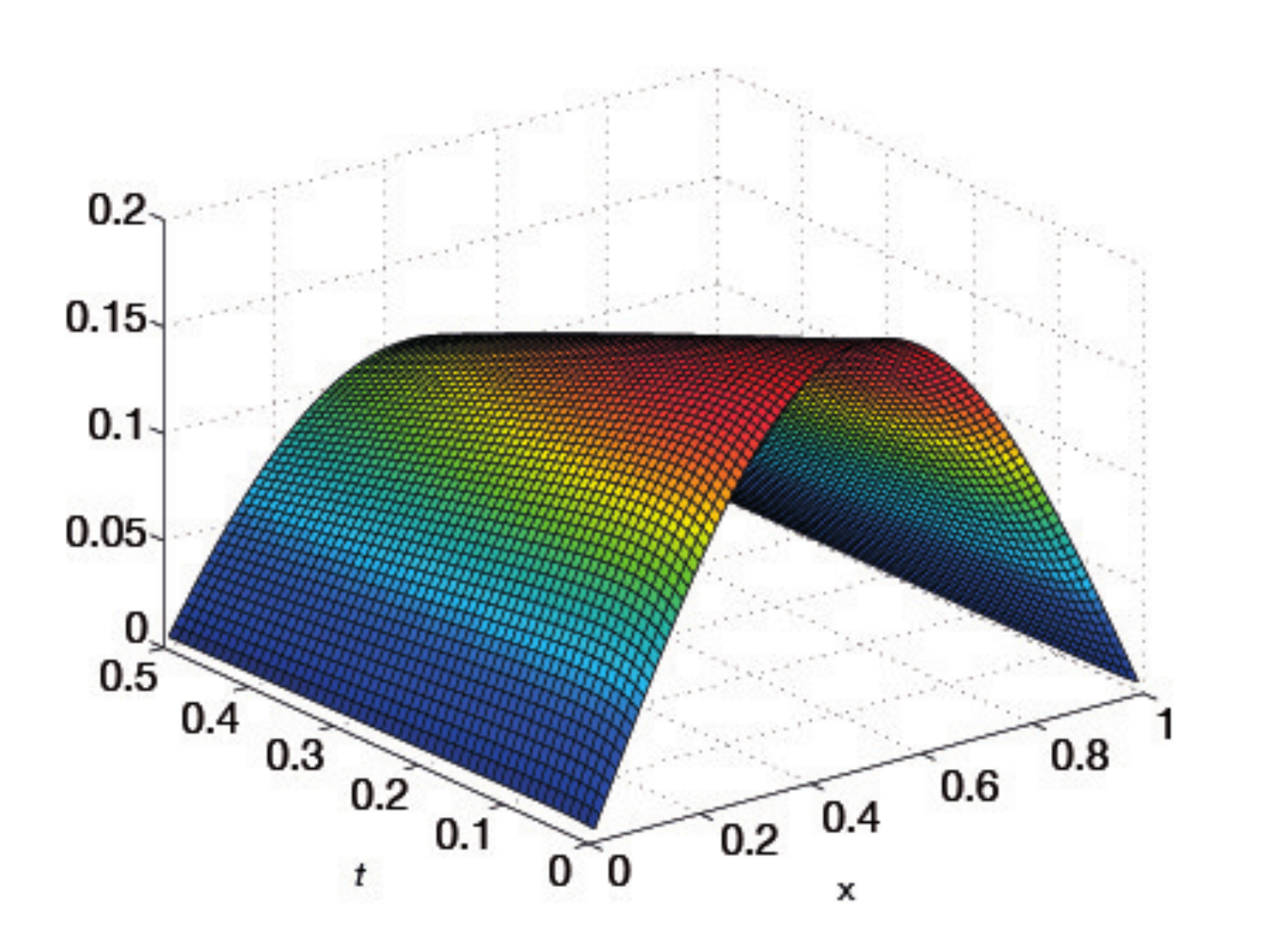}\hfill\includegraphics[height=43mm,width=40mm]{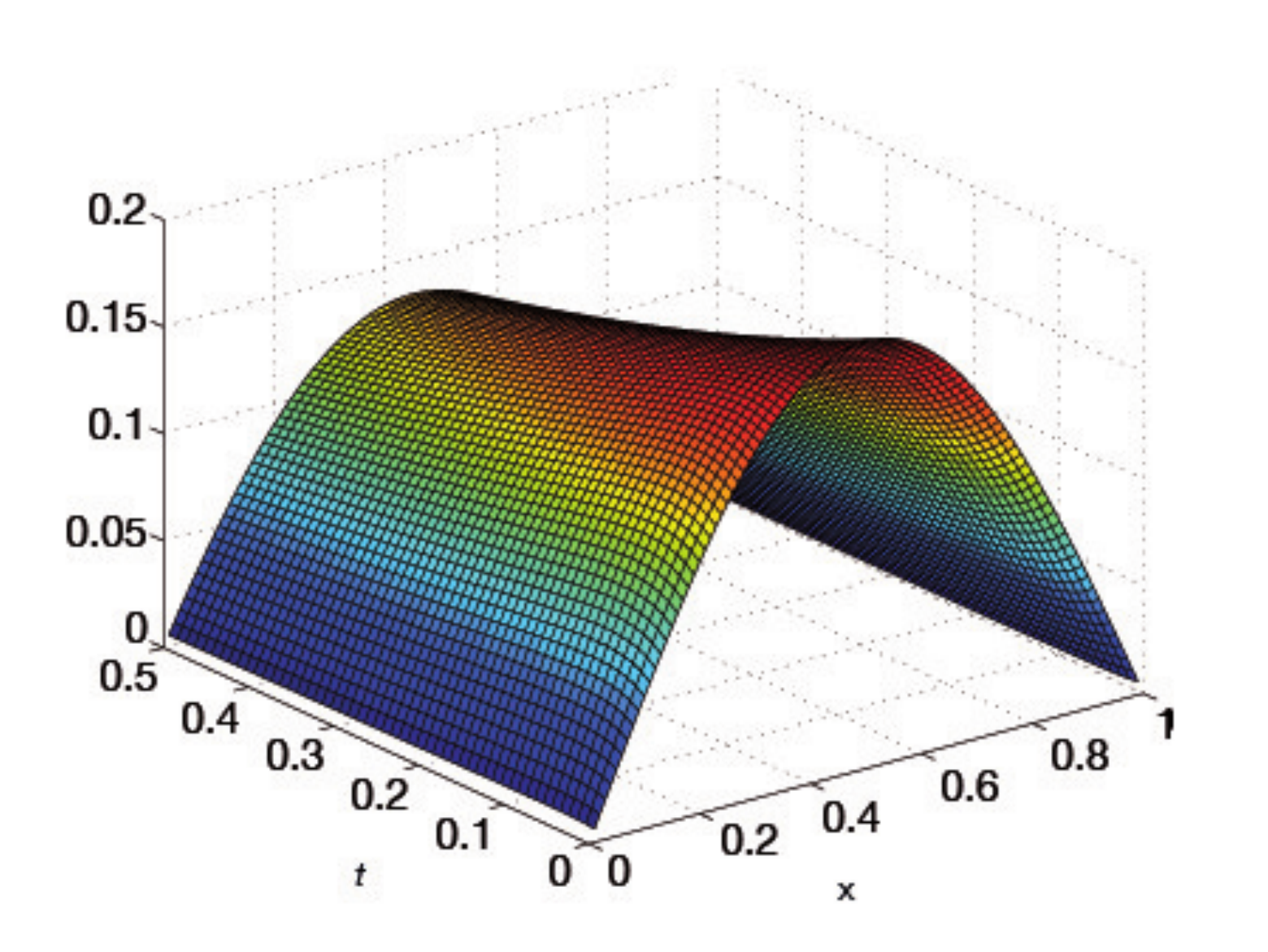}
\caption{Run~\ref{Run2}: NMPC state with $N=14$ (left plot), POD-NMPC state with $N=14$ (middle plot) and state with $u=-Ky$ (right plot)}
\label{pic:test2}
\end{figure}
The solution coming from the POD model is in the middle of Figure~\ref{pic:test2}. Note that $\mathcal E(\mathfrak l=3)=0.01, \mathcal E(\mathfrak l=13)=0,$ and $\mathsf{Err}(t;\mathfrak l)\le 10^{-3}$ for any $\mathfrak l$ and $t\ge t_\circ.$ Indeed, Table \ref{tab:test2} presents the evaluation of the cost functionals for the proposed algorithms and the CPU time which shows that the speed up by the reduced order approach is about 16.
\begin{table}[htbp]
\begin{center}
\begin{tabular}{lcccc}
\bottomrule
\hline                                                       
 & $\hat J$ & time & $K$ & $\|y^{FD}-y\|_{L^2(t_\circ,T;H)}$ \\
\midrule
Solution with $u=-Ky$ & 0.0035 &   & 1.50 & 0.0089\\
Alg.~\ref{Alg:NMPC} & 0.0027 & 65s &  & \\
Alg.~\ref{Alg:POD-NMPC} ($\mathfrak l=13$, $\mathfrak l^{DEIM}=15$) & 0.0032 & \phantom{1}5s &  & 0.0054\\
Alg.~\ref{Alg:POD-NMPC} ($\mathfrak l=3$, $\mathfrak l^{DEIM}=2$) & 0.0033 & \phantom{1}4s &  & 0.0055\\
\bottomrule
\end{tabular}
\end{center}
\caption{Run~\ref{Run2}: Evaluation of the cost functional, CPU times, suboptimal solution.}
\label{tab:test2}
\end{table}
Note that $K$ in Run~\ref{Run2} is smaller compared to Run~\ref{Run1} due to the constraint of the control space. Further, the error is presented in Table~\ref{tab:test2}. To study the influence of $\mathsf{Err}(t;\mathfrak l)$ we present in Figure~\ref{inf:test21}, on the left, how the optimal prediction horizon $N$ changes according to different tolerance. 
\begin{figure}[htbp]
\centering
\includegraphics[height=43mm,width=40mm]{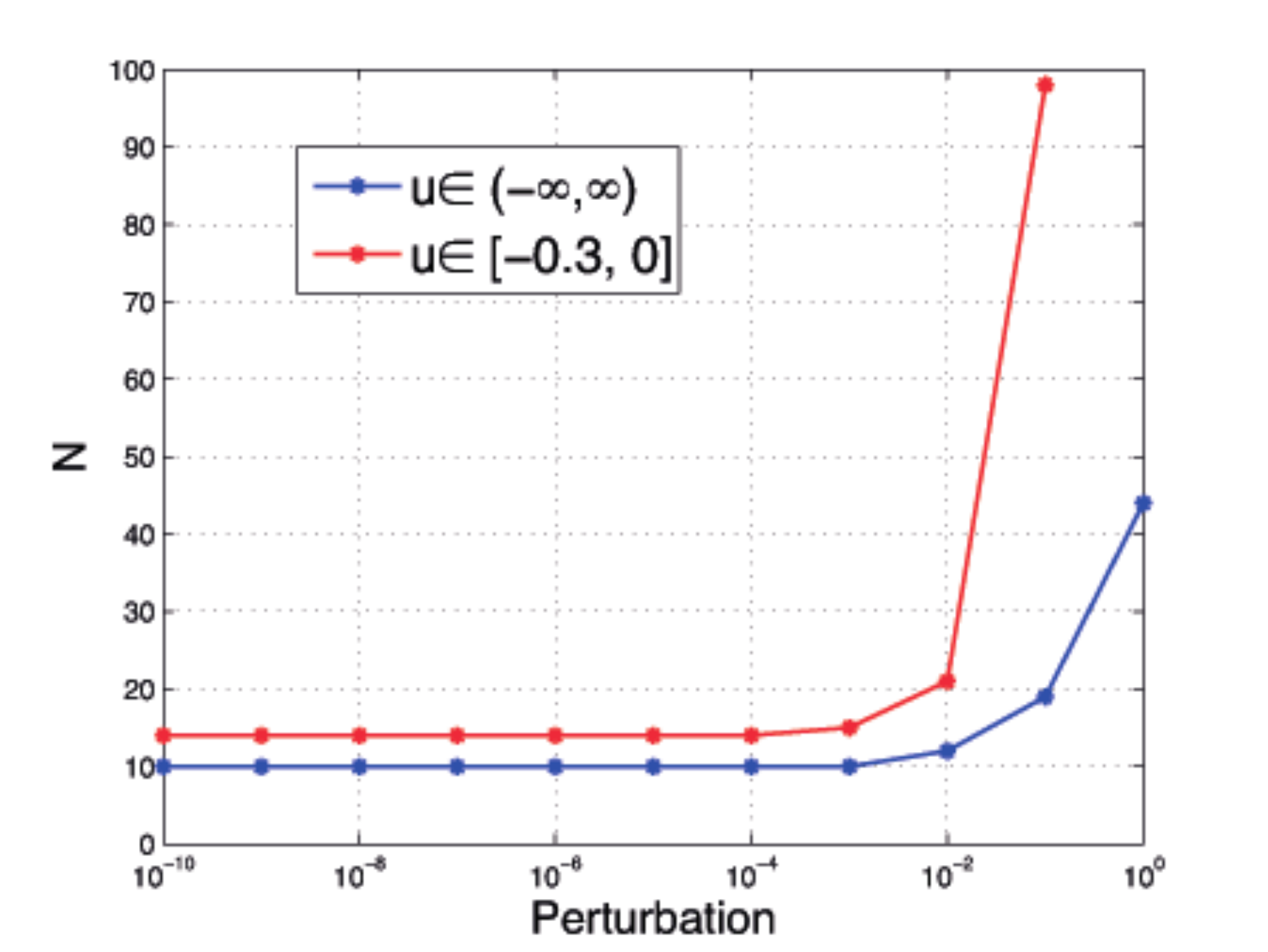}\hfill\includegraphics[height=43mm,width=40mm]{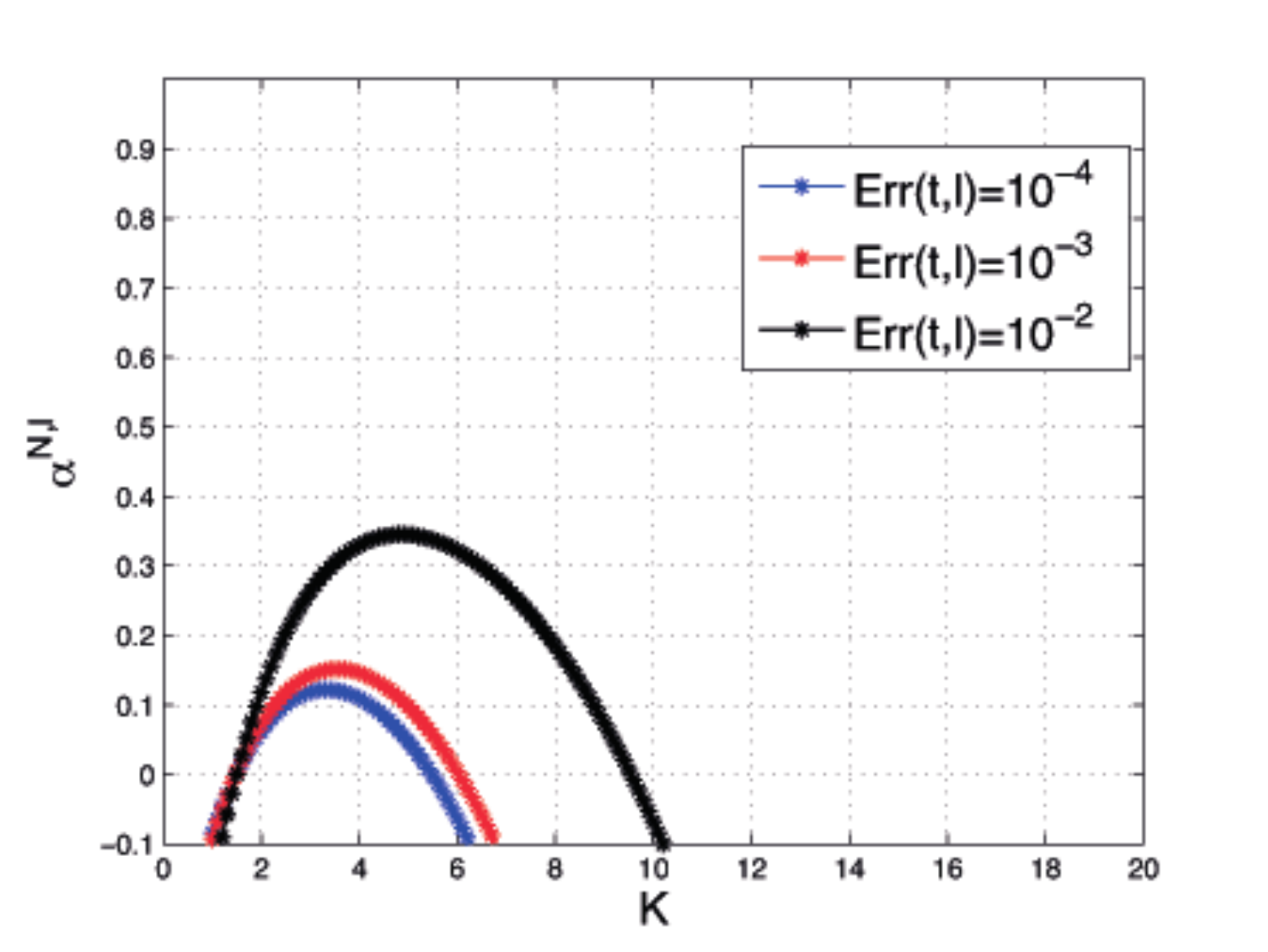}\hfill\includegraphics[height=43mm,width=40mm]{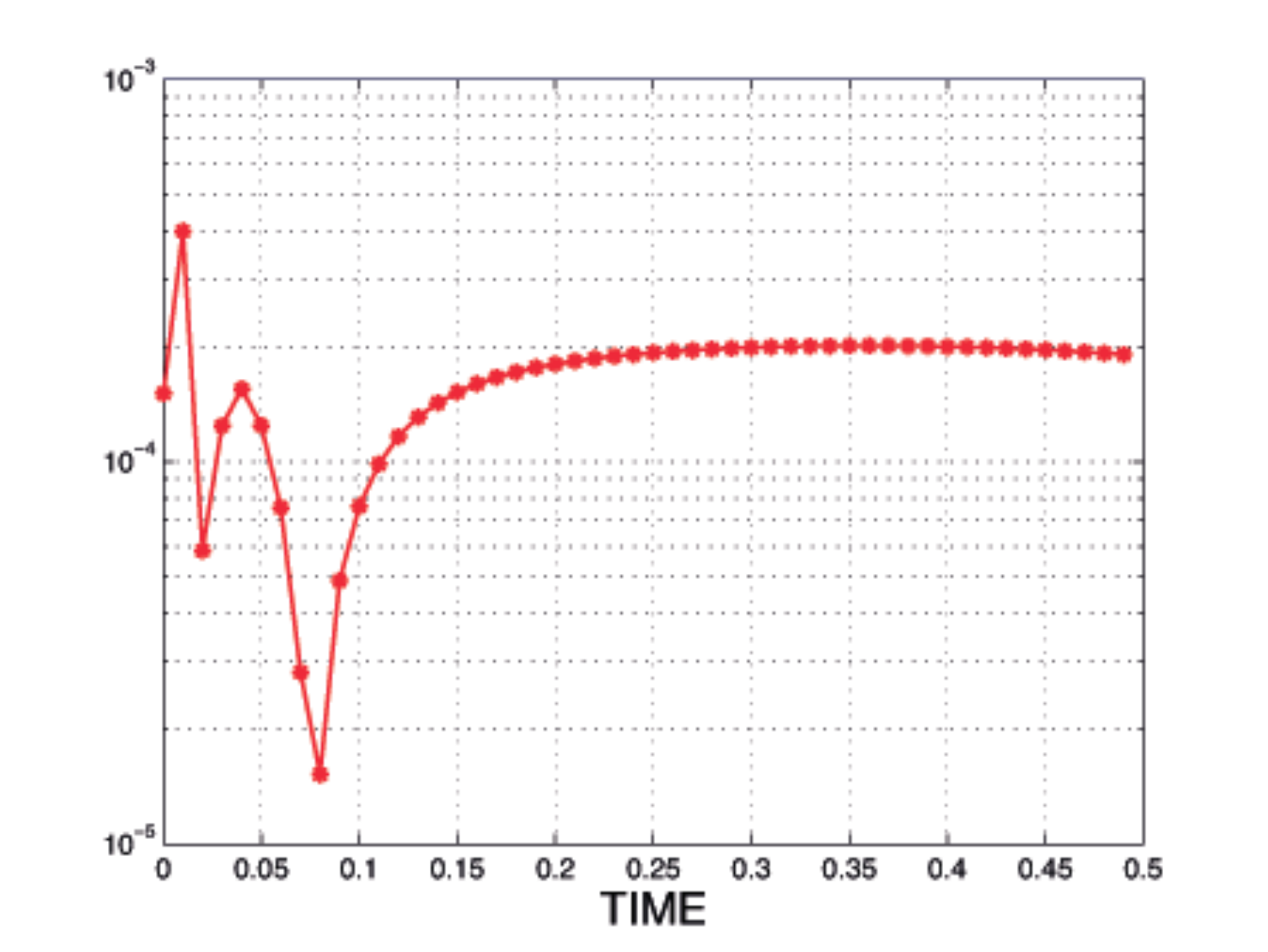}
\caption{Run~\ref{Run2}: Optimal horizon $N$ and $\alpha^{N,\mathfrak l}$ according to different $\mathsf{Err}(t;\mathfrak l)=10^{-3}$, Influence of the relative error $t\mapsto\mathsf{Err}(t;\mathfrak l)=10^{-3}$ for $\mathfrak l=3$.}
\label{inf:test21}
\end{figure}
The blue line corresponds to the optimal prediction horizon in Run~\ref{Run1}, and the red one to Run~\ref{Run2}. It turns out that, as long as $\mathsf{Err}(t;\mathfrak l)\le10^{-3},$ we can work exactly with the same horizon $N$ we had in the full model in both examples. In the middle plot of Figure~\ref{inf:test21} there is a zoom of the function $\alpha$ with different values of $\mathsf{Err}(t;\mathfrak l)$ with respect to Run~\ref{Run2}. The right plot of Figure~\ref{inf:test21} shows the relative error $\mathsf{Err}(t;\mathfrak l)$ for $0\le t\le 0.5$ with $\mathfrak l=3$. One of the big advantages of feedback control is the stabilization under perturbation of the system. The perturbation of the initial condition is a typical example which comes from many applications in fact, often the measurements may not be correct. For a given noise distribution $\delta=\delta(x)$ we consider a perturbation the following form:
\[
y_0(x)=\big(1+\delta(x)\big)y_\circ(x)\quad\text{for }x\in\Omega.
\]
The perturbation is applied only at every initial condition of the MPC algorithm (see \eqref{P-DynSys-FinHor} in Algorithm \ref{Alg:NMPC}) and it is random with respect to the spatial variable. The study of the asympotic stability does not change: we can compute the minimal prediction horizon as before. As we can see in Figure~\ref{fig:5} the POD-NMPC algorithm is able to stabilize with a noise of $|\delta (x)|\le30\%$.\hfill$\Diamond$
\begin{figure}[htbp]
\centering
\includegraphics[height=43mm,width=40mm]{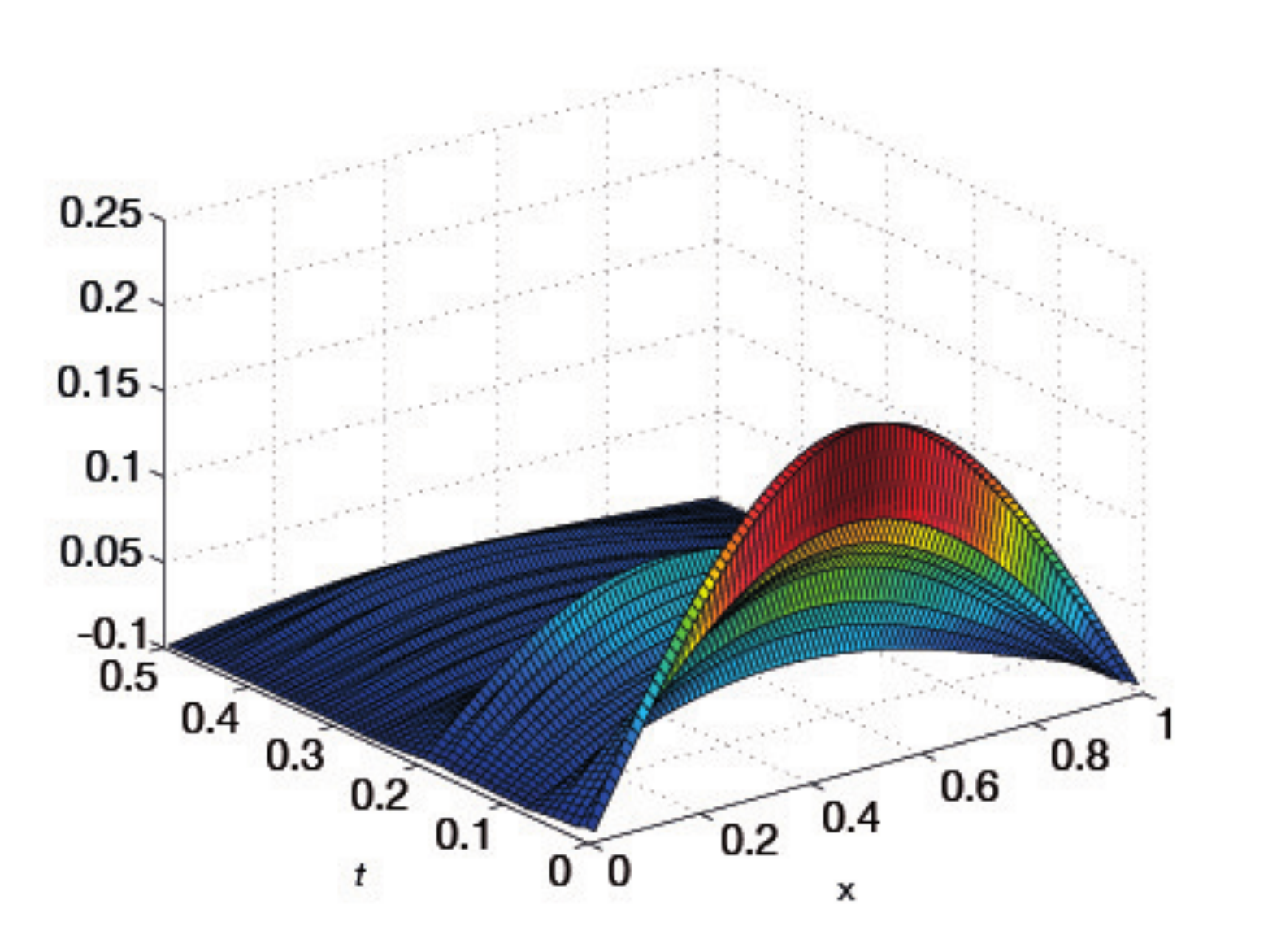}\hfill\includegraphics[height=43mm,width=40mm]{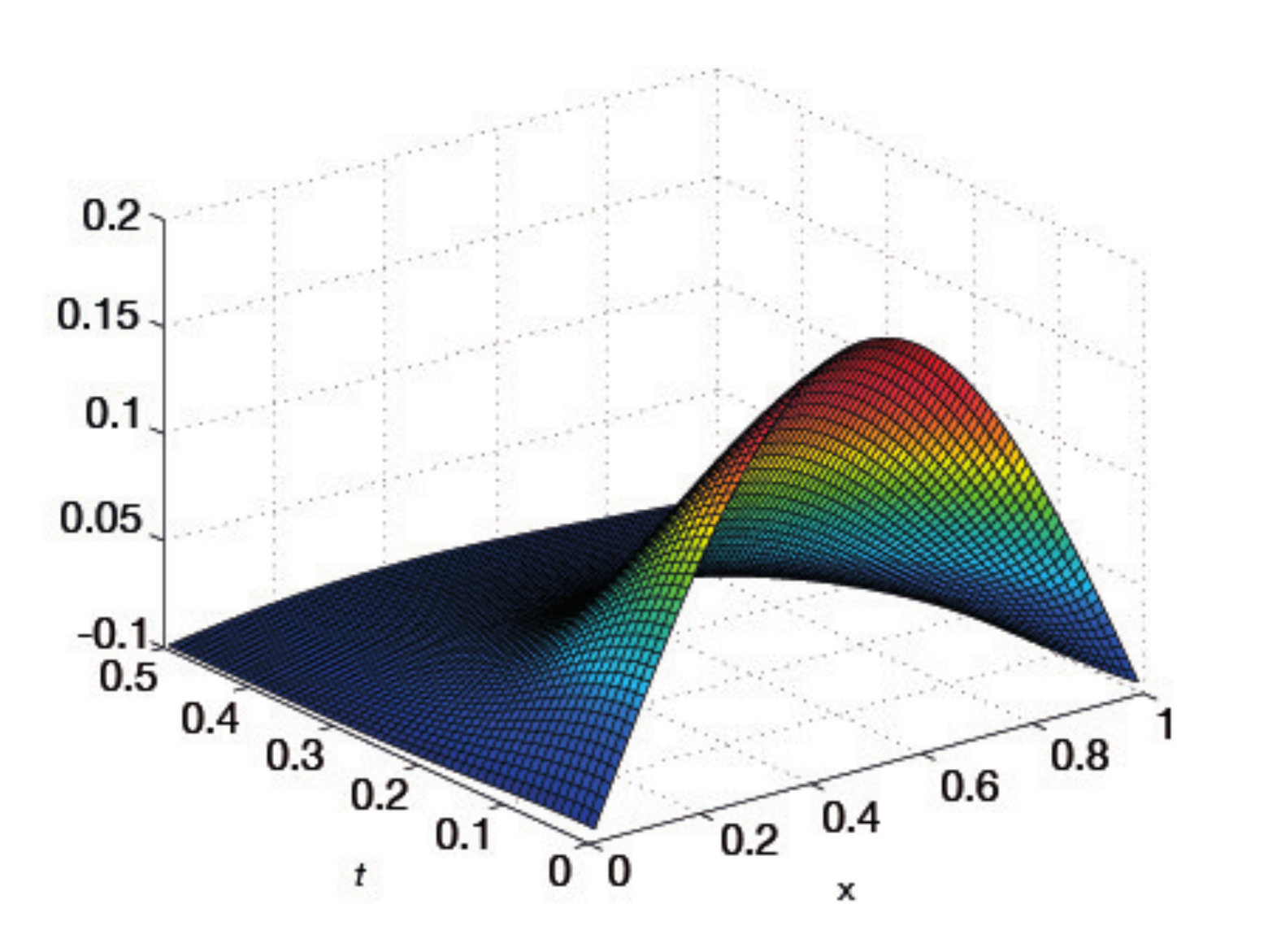}\hfill\includegraphics[height=43mm,width=40mm]{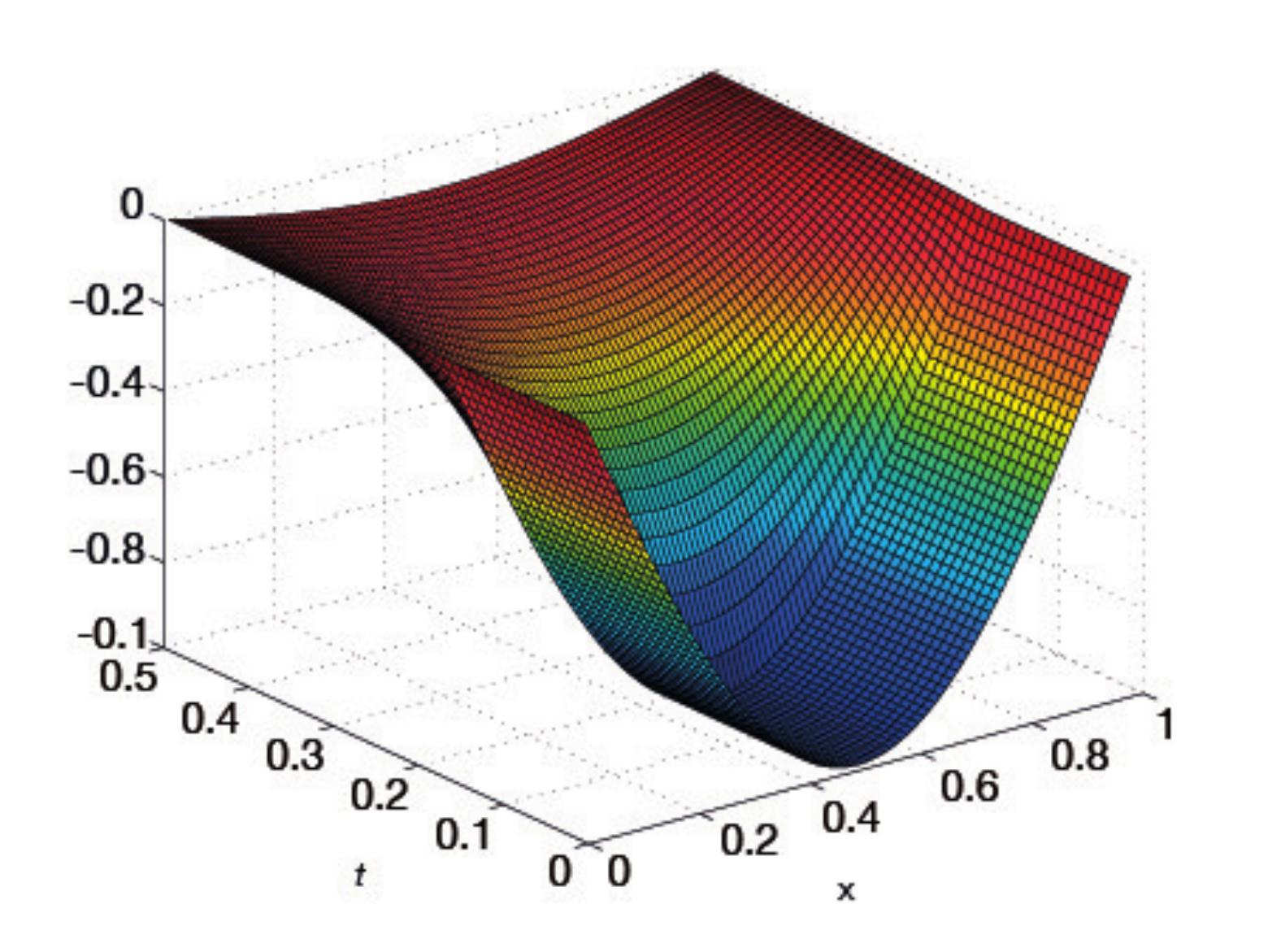}

\caption{Run~\ref{Run2}: POD-NMPC state with $30\%$ noise (left plot); Run~\ref{Run3}: NMPC state with $N=30$ (middle plot) and POD-NMPC state with $N=30$, $\mathfrak l=\mathfrak l^{DEIM}=16$ (right plot).}
\label{fig:5}
\end{figure}
\end{run}

\begin{run}[\em Constrained case with smooth initial data]
\label{Run3}
\rm
Now we decrease the diffusion term and, as a consequence, the prediction horizon $N$ increases; see Table \ref{par:3} and middle plot of Figure~\ref{fig:5}.
\begin{table}[htbp]
\begin{center}
\begin{tabular}{cccccccccc}
\toprule 
$T$ & $\Delta t$ & $\Delta x$ & $\theta$ & $\rho$ & $y_0(x)$ & $u_a$ & $u_b$ & $N$ & $K$\\
\midrule
0.5 & 0.01 & 0.01 & $1/\sqrt{2}$ & 10 & $0.2\sin(\pi x)$ & $-1$ & $0$ & 30& 5\\
\bottomrule
\end{tabular}
\end{center}
\caption{Run~\ref{Run3}: Setting for the optimal control problem.}
\label{par:3}
\end{table}
Even if the horizon is very large, the proposed Algorithm~\ref{Alg:POD-NMPC} accelerates the approximation of the problem. The decrease of $\theta$ may give some troubles with the POD-model since the domination of the convection term causes a high-variability in the solution, then a few basis functions will not suffice to obtain good surrogate models (see \cite{AF13,AF13bis}). Note that, in our example, the diffusion term is still relevant such that we can work with only 2 POD basis functions. The CPU time in the full model is $84$ seconds, whereas with a low-rank model, such as $\mathfrak l=2$ we obtained the solution in five seconds and an impressive speed up factor of  16. Even with a more accurate POD model we have a very good speed up factor of nine. The evaluation of the cost functional is given in Table \ref{tab:3}.
\begin{table}[htbp]
\begin{center}
\begin{tabular}{lcccc}
\toprule                                     
 & $\hat J$ & time & $K$ &$\|y^{FD}-y\|_{L^2(t_\circ,T;H)}$ \\
\midrule
Suboptimal solution ($u=-Ky$) & 0.0021 &   & 5 & 0.0208\\
Algorithm \ref{Alg:NMPC} & 0.0016 & 84s &  & \\
Algorithm \ref{Alg:POD-NMPC} ($\mathfrak l=16$, $\mathfrak l^{DEIM}=16$) & 0.0017 & \phantom{1}9s &  & 0.0092\\
Algorithm \ref{Alg:POD-NMPC} ($\mathfrak l=2$, $\mathfrak l^{DEIM}=3$) & 0.0018 & \phantom{1}5s &  & 0.0093\\
\bottomrule
\end{tabular}
\end{center}
\caption{Run~\ref{Run3}: Evaluation of the cost functional and CPU time. }
\label{tab:3}
\end{table}
In the right plot of Figure \ref{fig:5} the POD-NMPC state is plotted for $\mathfrak l=16$ POD basis and $\mathfrak l^{DEIM}=16$ DEIM ansatz functions. The error between the NMPC state and the POD-MPC state is less than 0.01 .\hfill$\Diamond$
\end{run}

\begin{run}[Constrained case with no-smooth initial data]
\label{Run4}
\rm
In the last test we focus on a different initial condition and different control constraints. The parameters are presented in Table \ref{par:4}.
\begin{table}[htbp]
\begin{center}
\begin{tabular}{cccccccccc}
\toprule
$T$ & $\Delta t$ & $\Delta x$ & $\theta$ & $\rho$ & $y_0(x)$ & $u_a$ & $u_b$ & $N$ & $K$\\
\midrule
0.5 & 0.01 & 0.01 & 1/2 & 5 & $0.1\sgn(x-0.3)$ & -1 & 1 & 43 & 9.99\\
\bottomrule
\end{tabular}
\end{center}
\caption{Run~\ref{Run4}: Setting for the optimal control problem.}
\label{par:4}
\end{table}
The minimal horizon $N$ which ensures asymptotic stability is $N=43$. Table \ref{tab:test4} emphazises again the performance of the POD-NMPC method with an acceleration 12 times faster than the full model.
\begin{table}[H]
\begin{center}
\begin{tabular}{lcccc}
\toprule                                                    
 & $\hat J$ & time & $K$ &$\|y^{FD}-y\|_{L^2(t_\circ,T;H)}$\\
\midrule
Solution with $u=-Ky$ & 4.7e-4 &   & 9.99& 0.0060\\
Alg.~\ref{Alg:NMPC} & 4.1e-4 & 50s  & \\
Alg.~\ref{Alg:POD-NMPC} ($\mathfrak l=17$, $\mathfrak l^{DEIM}=19$) & 4.4e-4 & 12s & & 0.0034\\
Alg.~\ref{Alg:POD-NMPC} ($\mathfrak l=3$, $\mathfrak l^{DEIM}=4$) & 4.4e-4 &  \phantom{1}4s&  &0.0035 \\
\bottomrule
\end{tabular}
\end{center}
\caption{Run~\ref{Run4}: Cost functional, CPU time and suboptimal solution.}
\label{tab:test4}
\end{table}
The evaluation of the cost functional gives the same order in all the simulation we provide. In Figure \ref{fig:4} we present the NMPC state for $N=43$ (left plot), the POD-NMPC state with $N=43$, $\mathfrak l=3$, $\mathfrak l^{DEIM}=4$ (middle plot) and the increase of the optimal horizon $N$ according to the perturbation $\mathsf{Err}(t;\mathfrak l)$.
\begin{figure}[htbp]
\centering
\includegraphics[height=43mm,width=37mm,clip=true]{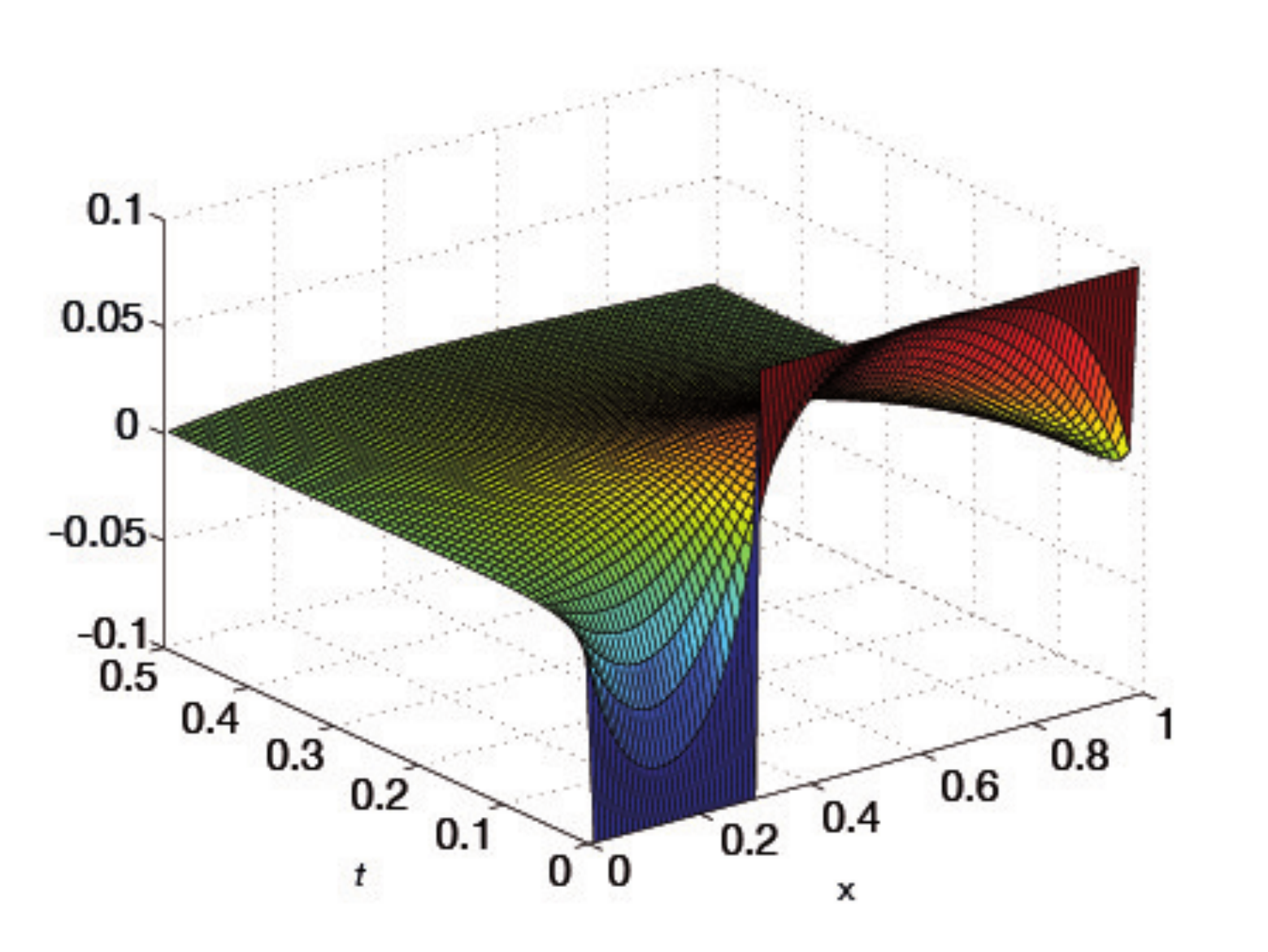}\hfill\includegraphics[clip=true,height=43mm,width=37mm]{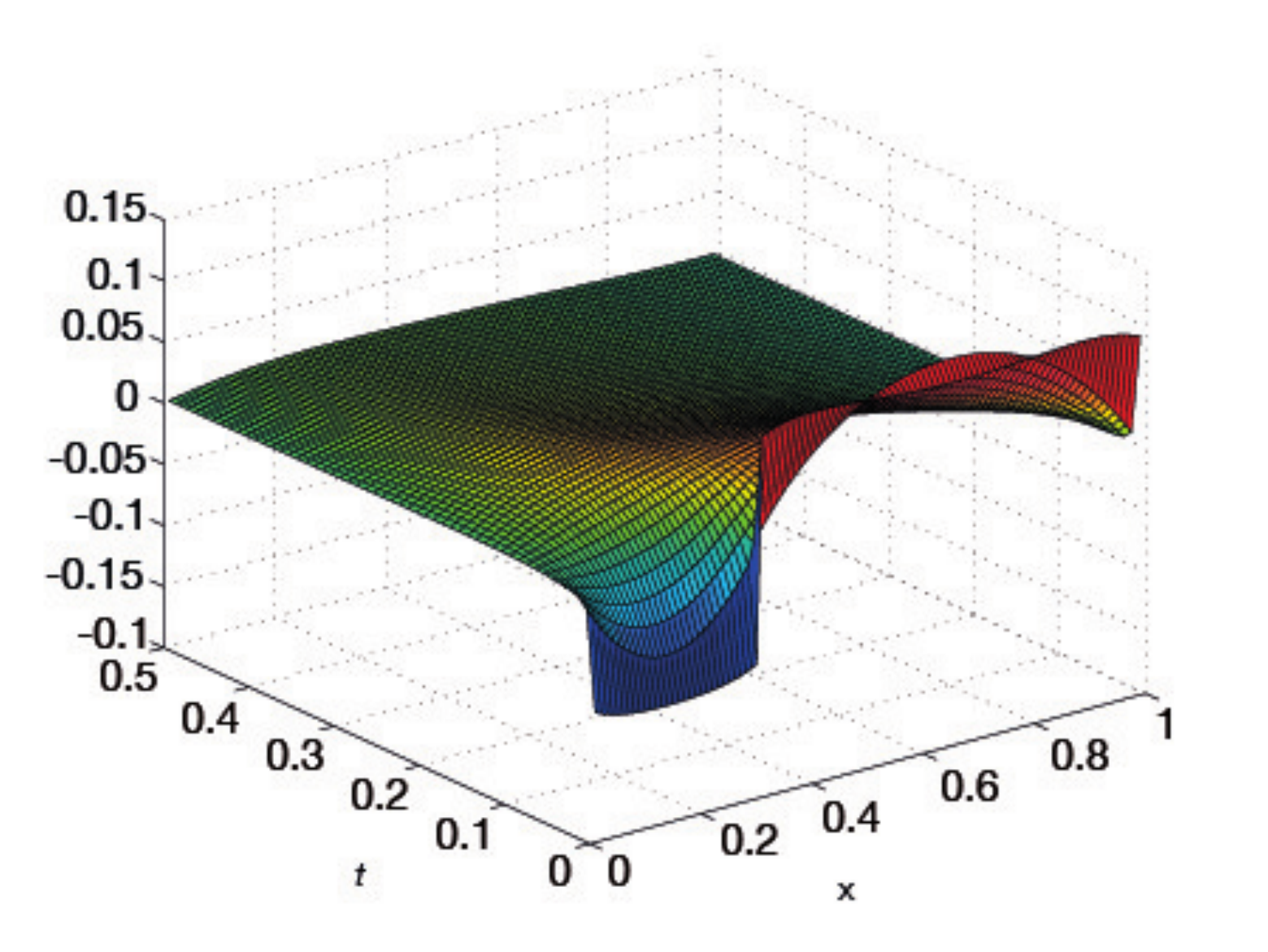}\hfill\includegraphics[height=43mm,width=37mm]{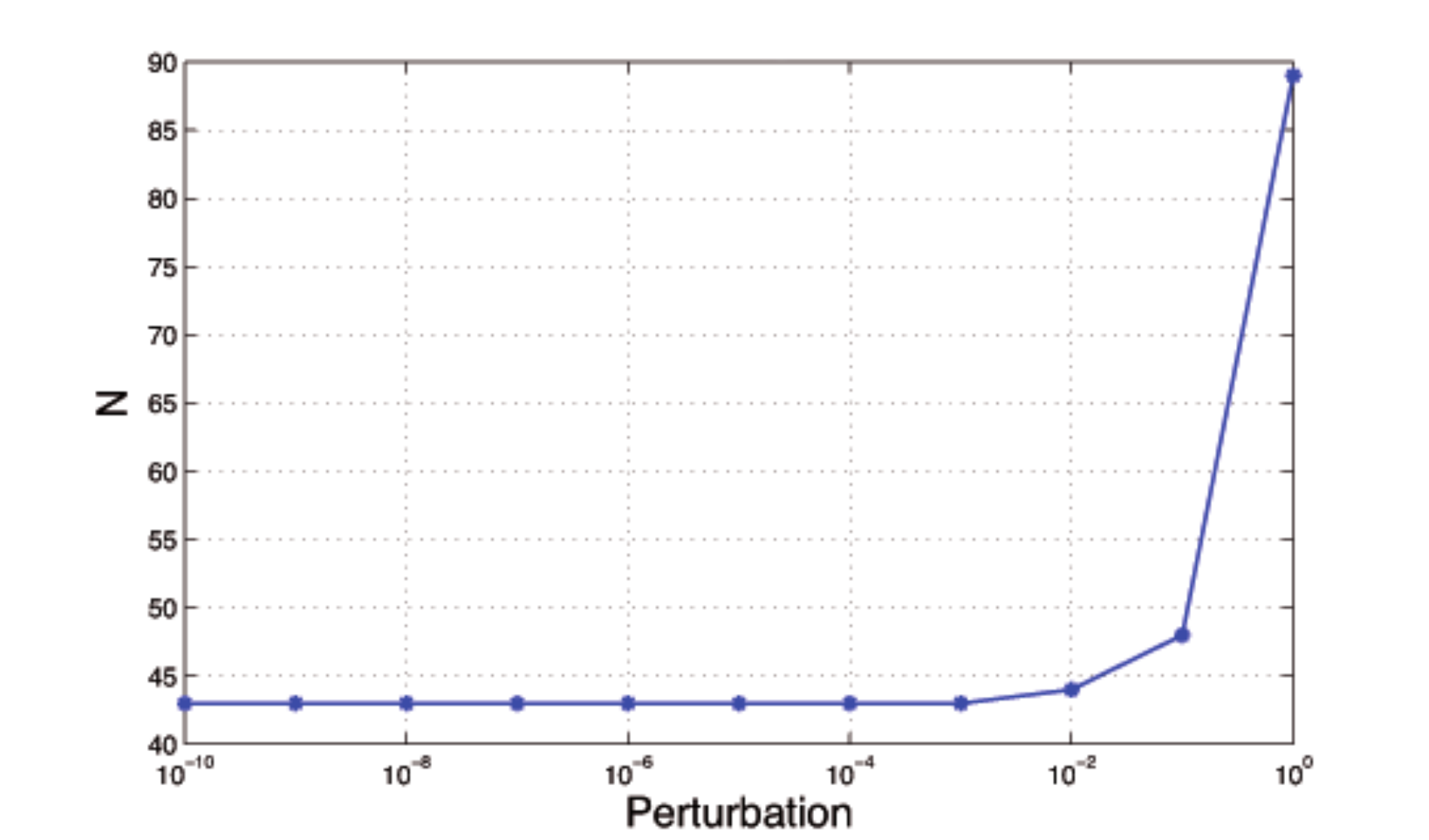}

\caption{Run~\ref{Run4}: NMPC state with $N=43$ (left plot), POD-NMPC state with $N=43,$ $\mathfrak l=17$, $\mathfrak l^{DEIM}=19$ (middle plot) and increase of the optimal horizon $N$ according to the perturbation $\mathsf{Err}(t;\mathfrak l)$ (right plot).} 
\label{fig:4}
\end{figure}
The error between the NMPC state and the POD-MPC state is 0.0035 when $\mathcal E(\mathfrak l=3)=0.01,$ whereas for $\mathcal E(\mathfrak l=17)=0$ the error is 0.0034. \hfill$\Diamond$
\end{run}

\section{Conclusions}

We have proposed a new numerical method for optimal control problems which tries to stabilize a one dimensional semilinear parabolic equation by means of Nonlinear MPC. We presented asymptotic stability conditions, where the control space is bounded for a suboptimal problem coming from a particular class of feedback controls.\\
Since the CPU time of the full dimensional algorithm may increase with the dimension of the prediction horizon, we have presented a deep study of the suboptimal model which comes from POD model reduction. We have given an {\em a-priori} error estimate for the computation of the prediction horizon of the suboptimal model. The new reduced model approach turns out to be computationally very efficient with respect to the full dimensional problem. If the approximation quality \eqref{Eq:ErrTerm} of the reduced-order model is taken into account, stabilization is also guaranteed by our theory. Although the algorithm is applied to a one dimensional problem, the theory is rather general and can be applied to higher dimensional equations, not only with POD model reduction but any (reduced-order) method provided the error term in \eqref{Eq:ErrTerm} is small for reasonable small $\mathfrak l$.


\begin{appendix}

\section{Proof of Lemma~\ref{lem:stab}}

\label{App-Lemma3.6}

Choosing $u(t)=-Ky(t)$, $\varphi=y(t)$ in \eqref{PDEWeakForm} and using $\int_\Omega y(t)^4\,\mathrm dx\ge 0$, $\int_\Omega y_x(t)y(t)\,\mathrm dx=0$ f.a.a. $t\ge t_\circ$ we find
\begin{align*}
\frac{1}{2}\frac{\mathrm d}{\mathrm dt}\,{\| y(t)\|}^2_H+\theta\,{\|y(t)\|}_V^2+(K-\rho)\,{\|y(t)\|}_H^2\le0\quad\text{f.a.a. }t\ge t_\circ
\end{align*}
Hence, \eqref{Eq-ocp.embed} imply
\[
\frac{\mathrm d}{\mathrm dt}\,{\| y(t)\|}^2_H\le-2\bigg(\frac{\theta}{C_V}+K-\rho\bigg)\,{\|y(t)\|}_H^2=-2\gamma(K)\,{\|y(t)\|}_H^2\quad\text{f.a.a. }t\ge t_\circ.
\]
Thus, by Gronwall's inequality we derive
\[
{\| y(t)\|}^2_H\le e^{-2\gamma(K)(t-t_\circ)}\,{\|y_\circ\|}_H^2\quad\text{f.a.a. }t\ge t_\circ.
\]
which gives \eqref{AprioriEst}.\hfill$\Box$


\section{Proof of Theorem~\ref{Pro:Apriori}}

\label{AppA}

\noindent
Recall that $H^\mathfrak l\subset V$ holds. Consequently, $\|\psi_i^H-\mathcal P^\mathfrak l_H\psi_i^H\|_V$ is well-defined for $1\le i\le\mathfrak l$. First we review a result from \cite[Theorem~6.2]{Sin13}, which is essential in our proof of Theorem~\ref{Pro:Apriori} for the choice $X=H$: Suppose that $y^k\in L^2(t_\circ,t_\circ^N;V)$ for $1\le k\le \wp$. Then,
\begin{equation}
\label{GVLuminy:Theorem3.2.1}
\sum_{k=1}^\wp\int_{t_\circ}^{t_\circ^N} {\|y^k(t)-\mathcal P_H^\mathfrak l y^k(t)\|}_V^2\,\mathrm dt=\sum_{i=\mathfrak l+1}^d\lambda_i^H\,{\|\psi_i^H-\mathcal P^\mathfrak l_H\psi_i^H\|}_V^2.
\end{equation}
Moreover, $\mathcal P^\mathfrak l_H y^k$ converges to $y^k$ in $L^2(0,T;V)$ as $\mathfrak l$ tends to $\infty$ for each $k\in\{1,\ldots,\wp\}$.

\smallskip\noindent{\em Proof of Theorem~\ref{Pro:Apriori}.} To derive an error estimate for $\|y-y^\mathfrak l\|_{\mathbb Y^N(t_\circ)}$ we
make use of the decomposition
\[
y(t)-y^\mathfrak l(t)=y(t)-\mathcal P^\mathfrak ly(t)+\mathcal P^\mathfrak l y(t)-y^\mathfrak l(t)=\varrho^\mathfrak l(t)+\vartheta^\mathfrak l(t)\text{ f.a.a. }t\in[t_\circ,t_\circ^N]
\]
with $\varrho^\mathfrak l(t)=y(t)-\mathcal P^\mathfrak ly(t)\in(X^\mathfrak l)^\bot$ and $\vartheta^\mathfrak l(t)=\mathcal P^\mathfrak l y(t)-y^\mathfrak l(t)\in X^\mathfrak l$. Recall that $\mathbb Y^N(t_\circ)=W(t_\circ,t_\circ^N)$ holds. Since $\varrho_t^\mathfrak l(t)\in V$ holds f.a.a. $t\in[t_\circ,t_\circ^N]$, we have $\|\varrho_t^\mathfrak l(t)\|_{V'}=\|\varrho_t^\mathfrak l(t)\|_V$ due to the Riesz theorem \cite[p.~43]{RS80}. Hence it follows from \eqref{Eq:RhoError-1} and \eqref{GVLuminy:Theorem3.2.1} that
\begin{equation}
\label{Eq:RhoError-2}
\begin{aligned}
{\|\varrho^\mathfrak l\|}_{\mathbb Y^N(t_\circ)}^2&=\int_{t_\circ}^{t_\circ^N}{\|\varrho^\mathfrak l(t)\|}_V^2+{\|\varrho^\mathfrak l_t(t)\|}_V^2\,\mathrm dt\\
&=\left\{
\begin{aligned}
&\sum_{i=\mathfrak l+1}^d\lambda_i^V&&\text{for }X=V,\\
&\sum_{i=\mathfrak l+1}^d\lambda_i^H{\|\psi_i^H-\mathcal P_H^\mathfrak l\psi_i^H\|}_V^2&&\text{for }X=H.
\end{aligned}
\right.
\end{aligned}
\end{equation}
Next we estimate $\vartheta^\mathfrak l(t)$. We infer from $\vartheta^\mathfrak l(t)=\mathcal P^\mathfrak l y(t)-y^\mathfrak l(t)$ that
\begin{equation}
\label{App-1}
\begin{aligned}
&{\langle \vartheta^\mathfrak l_t(t),\psi\rangle}_{V',V}+{\langle\theta \vartheta^\mathfrak l(t),\psi\rangle}_V\\
&\quad={\langle y_t(t)-y^\mathfrak l_t(t)+\mathcal P^\mathfrak ly(t)-y^\mathfrak l(t),\psi\rangle}_{V',V}+\theta\,{\langle \mathcal P^\mathfrak ly(t)-y^\mathfrak l(t),\psi\rangle}_V
\end{aligned}
\end{equation}
for all $\psi\in X^\mathfrak l$ and f.a.a. $t\in[t_\circ,t_\circ^N]$. For $X=V$ we have 
\[
{\langle\mathcal P^\mathfrak l_Vy(t),\psi\rangle}_V={\langle y(t),\psi\rangle}_V\quad\text{for all }\psi\in X^\mathfrak l\text{ and f.a.a. }t\in[t_\circ,t_\circ^N].
\]
Hence, we derive from \eqref{App-1}, \eqref{StateEquation} and \eqref{Eq:PODGalState} that
\begin{equation}
\label{App-2}
\begin{aligned}
&{\langle \vartheta^\mathfrak l_t(t),\psi\rangle}_{V',V}+{\langle\theta \vartheta^\mathfrak l(t),\psi\rangle}_V\\
&\quad={\langle \rho(y(t)-y^\mathfrak l(t))-\rho(y(t)^3-y^\mathfrak l(t)^3),\psi\rangle}_H+{\langle \mathcal P^\mathfrak l_V y_t(t)-y_t(t),\psi\rangle}_{V',V}
\end{aligned}
\end{equation}
for all $\psi\in V^\mathfrak l$ and f.a.a. $t\in[t_\circ,t_\circ^N]$. For $s\in [0,1]$ we define the function $\xi^\mathfrak l(s)=y^\mathfrak l+s(y-y^\mathfrak l)$. Then it follows from \eqref{APrioriState} and \eqref{Eq:LinfEstPOD} that
\[
{\|\xi^\mathfrak l(s)\|}_{L^\infty(Q^N)}\le s\,{\|y\|}_{L^\infty(Q^N)}+(1-s)\,{\|y^\mathfrak l\|}_{L^\infty(Q^N)}\le C_1\quad\text{for all }s\in[0,1]
\]
with a constant $C_1>0$ dependent on $y_\circ$, $u_a$ and $u_b$, but independent of $y$, $y^\mathfrak l$ and $\mathfrak l$. By the mean value theorem we obtain
\begin{align*}
{\langle y(t)^3-y^\mathfrak l(t)^3,\psi\rangle}_H&=\bigg\langle\frac{1}{4}\int_0^1\xi^\mathfrak l(s;t)^4\big(y(t)-y^\mathfrak l(t)\big)\,\mathrm ds,\psi\bigg\rangle_H\\
&\le C_2\,{\|y(t)-y^\mathfrak l(t)\|}_H{\|\psi\|}_H\quad\text{for all } \psi\in V^\mathfrak l
\end{align*}
with $C_2=C_1^4/4$. We set $C_3=\rho(1+C_2)$. Hence, choosing $\psi=\vartheta^\mathfrak l(t)\in V^\mathfrak l$ and utilizing $\theta\ge\theta_a>0$ we obtain from $y(t)-y^\mathfrak l(t)=\varrho^\mathfrak l(t)+\vartheta^\mathfrak l(t)$, \eqref{App-2}, \eqref{Eq-ocp.embed} and Young's inequality
\begin{align*}
&\frac{1}{2}\frac{\mathrm d}{\mathrm dt}\,{\|\vartheta^\mathfrak l(t)\|}_H^2+\theta_a\,{\|\vartheta^\mathfrak l(t)\|}_V^2\\
&\le C_3\,{\|y(t)-y^\mathfrak l(t)\|}_H{\|\vartheta^\mathfrak l(t)\|}_H+{\|\mathcal P^\mathfrak l_V y_t(t)-y_t(t)\|}_{V'}{\|\vartheta^\mathfrak l(t)\|}_V\\
&\le C_3\,\bigg(\frac{C_V^2}{2}\,{\|\varrho^\mathfrak l(t)\|}_V^2+\frac{3}{2}\,{\|\vartheta^\mathfrak l(t)\|}_H^2\bigg)+\frac{1}{2\theta_a}{\|\varrho^\mathfrak l_t(t)\|}_{V'}^2+\frac{\theta_a}{2}\,{\|\vartheta^\mathfrak l(t)\|}_V^2\\
&\le\frac{C_4}{2}\,\big({\|\varrho^\mathfrak l(t)\|}_V^2+{\|\varrho^\mathfrak l_t(t)\|}_V^2\big)+\frac{3C_3}{2}\,{\|\vartheta^\mathfrak l(t)\|}_H^2+\frac{\theta_a}{2}\,{\|\vartheta^\mathfrak l(t)\|}_V^2
\end{align*}
f.a.a. $t\in[t_\circ,t_\circ^N]$ with the constant $C_4=\max(C_3C_V^2,1/\theta_a)$. Hence, we have
\begin{equation}
\label{App-3}
\frac{\mathrm d}{\mathrm dt}\,{\|\vartheta^\mathfrak l(t)\|}_H^2+\theta_a\,{\|\vartheta^\mathfrak l(t)\|}_V^2\le C_5\,\big({\|\varrho^\mathfrak l(t)\|}_V^2+{\|\varrho^\mathfrak l_t(t)\|}_{V}^2+{\|\vartheta^\mathfrak l(t)\|}_H^2\big)
\end{equation}
for $C_5=\max(C_4,3C_3)$ and f.a.a. $t\in[t_\circ,t_\circ^N]$. By Gronwall's inequality and \eqref{Eq:RhoError-2} we derive from \eqref{App-3}
\begin{equation}
\label{App-4}
\begin{aligned}
{\|\vartheta^\mathfrak l(t)\|}_H^2&\le e^{C_5(t-t_\circ)}\bigg({\|\vartheta^\mathfrak l(t_\circ)\|}_H^2+\int_{t_\circ}^{t_\circ^N}{\|\varrho^\mathfrak l(s)\|}_V^2+{\|\varrho^\mathfrak l_t(s)\|}_V^2\,\mathrm ds\bigg)\\
&\le C_6\bigg({\|\vartheta^\mathfrak l(t_\circ)\|}_H^2+\sum_{i=\mathfrak l+1}^\infty\lambda_i^V\bigg)\quad\text{f.a.a. }t\in[t_\circ,t_\circ^N]
\end{aligned}
\end{equation}
f.a.a. $t\in[t_\circ,t_\circ^N]$ with $C_6=e^{C_5(t_\circ^N-t_\circ)}$. Now we turn to the case $X=H$. We have
\[
{\langle\mathcal P^\mathfrak l_Hy(t),\psi\rangle}_V={\langle y(t),\psi\rangle}_V+{\langle\mathcal P^\mathfrak l_Hy(t)-y(t),\psi\rangle}_V
\]
so that \eqref{App-1}, \eqref{StateEquation} and \eqref{Eq:PODGalState} that
\begin{equation}
\label{App-5}
\begin{aligned}&{\langle \vartheta^\mathfrak l_t(t),\psi\rangle}_{V',V}+{\langle\theta \vartheta^\mathfrak l(t),\psi\rangle}_V={\langle \rho(y(t)-y^\mathfrak l(t))-\rho(y(t)^3-y^\mathfrak l(t)^3),\psi\rangle}_H\\
&\hspace{30mm}+{\langle \mathcal P^\mathfrak l_H y(t)-y(t),\psi\rangle}_V+{\langle \mathcal P^\mathfrak l_V y_t(t)-y_t(t),\psi\rangle}_{V',V}
\end{aligned}
\end{equation}
for all $\psi\in H^\mathfrak l$ and f.a.a. $t\in[t_\circ,t_\circ^N]$. Now we proceed analogously as in the case $X=V$ and obtain
\begin{align*}
&\frac{1}{2}\frac{\mathrm d}{\mathrm dt}\,{\|\vartheta^\mathfrak l(t)\|}_H^2+\theta_a\,{\|\vartheta^\mathfrak l(t)\|}_V^2\\
&\le\frac{C_7}{2}\,\big({\|\varrho^\mathfrak l(t)\|}_V^2+{\|\varrho^\mathfrak l_t(t)\|}_V^2+{\|\vartheta^\mathfrak l(t)\|}_H^2\big)+\frac{\theta_a}{2}\,{\|\vartheta^\mathfrak l(t)\|}_V^2
\end{align*}
f.a.a. $t\in[t_\circ,t_\circ^N]$ with the constant $C_7=\max(C_3C_V^2+2/\theta_a,2C_3)$. Therefore, we derive 
\begin{equation}
\label{App-6}
\frac{\mathrm d}{\mathrm dt}\,{\|\vartheta^\mathfrak l(t)\|}_H^2+\theta_a\,{\|\vartheta^\mathfrak l(t)\|}_V^2\le C_7\,\big({\|\varrho^\mathfrak l(t)\|}_V^2+{\|\varrho^\mathfrak l_t(t)\|}_V^2+{\|\vartheta^\mathfrak l(t)\|}_H^2\big);
\end{equation}
compare \eqref{App-3}. Utilizing Gronwall's inequality and \eqref{Eq:RhoError-2} we infer -- instead of \eqref{App-4} -- that
\begin{equation}
\label{App-7}
{\|\vartheta^\mathfrak l(t)\|}_H^2\le C_8\bigg({\|\vartheta^\mathfrak l(t_\circ)\|}_H^2+\sum_{i=\mathfrak l+1}^\infty\lambda_i^H\,{\|\psi_i^H-\mathcal P_H^\mathfrak l\psi_i^H\|}_V^2\bigg)
\end{equation}
f.a.a. $t\in[t_\circ,t_\circ^N]$ with $C_8=e^{C_7(t_\circ^N-t_\circ)}$.  We summarize \eqref{App-4} and \eqref{App-6} in
\begin{equation}
\label{Eq:HEst}
{\|\vartheta^\mathfrak l(t)\|}_H^2\le C_9\cdot\left\{
\begin{aligned}
&{\|\vartheta^\mathfrak l(t_\circ)\|}_H^2+\sum_{i=\mathfrak l+1}^d\lambda_i^V&&\text{for }X=V,\\
&{\|\vartheta^\mathfrak l(t_\circ)\|}_H^2+\sum_{i=\mathfrak l+1}^d\lambda_i^H{\|\psi_i^H-\mathcal P_H^\mathfrak l\psi_i^H\|}_V^2&&\text{for }X=H
\end{aligned}
\right.
\end{equation}
f.a.a. $t\in[t_\circ,t_\circ^N]$ with $C_9=\max(C_6,C_8)$. Furthermore, \eqref{App-3} and \eqref{App-6}, respectively, imply by integration over $[t_\circ,t_\circ^N]$
\begin{equation}
\label{Eq:thetaL2V}
\begin{aligned}
&{\|\vartheta^\mathfrak l\|}_{L^2(t_\circ,t_\circ^N;V)}^2\le \frac{1}{\theta_a}\,{\|\vartheta^\mathfrak l
(t_\circ)\|}_H^2+\frac{C_{10}}{\theta_a}\,\big({\|\varrho^\mathfrak l\|}_{W(t_\circ,t_\circ^N)}^2+{\|\vartheta^\mathfrak l\|}_{L^2(t_\circ,t_\circ^N;H)}^2\big)\\
&\quad\le C_{11}\cdot\left\{
\begin{aligned}
&{\|\vartheta^\mathfrak l(t_\circ)\|}_H^2+\sum_{i=\mathfrak l+1}^d\lambda_i^V&&\text{for }X=V,\\
&{\|\vartheta^\mathfrak l(t_\circ)\|}_H^2+\sum_{i=\mathfrak l+1}^d\lambda_i^H{\|\psi_i^H-\mathcal P_H^\mathfrak l\psi_i^H\|}_V^2&&\text{for }X=H
\end{aligned}
\right.
\end{aligned}
\end{equation}
with $C_{10}=\max(C_5,C_7)$ and $C_{11}=C_{10}\max(1,(t_\circ^N-t_\circ)\max(C_6,C_8))/\theta_a$. From estimates \eqref{Eq:HEst}, \eqref{Eq:thetaL2V}, from
\[
y(t)^3-y^\mathfrak l(t)^3=\big(\varrho^\mathfrak l(t)+\vartheta^\mathfrak l(t)\big)\big(y(t)^2+y(t)y^\mathfrak l(t)+y^\mathfrak l(t)^2\big)
\quad \text{f.a.a. }t\in[t_\circ,t_\circ^N]
\]
and from the embedding inequalities \cite{Eva08}
\begin{align}
\nonumber
{\|\varphi\|}_{L^\infty(\Omega)}&\le C_\infty\,{\|\varphi\|}_V&&
\text{for all }\varphi\in V,\\
\label{Embed-2}
{\|\varphi\|}_{C([t_\circ,t_\circ^N];H)}&\le C_W\,{\|\varphi\|}_{W(t_\circ,t_\circ^N)}&&
\text{for all }\varphi\in W(t_\circ,t_\circ^N)
\end{align}
for two constants $C_\infty,\,C_W>0$ we infer that
\begin{align*}
&{\|\vartheta^\mathfrak l\|}_{L^2(t_\circ,t_\circ^N;V')}=\sup_{{\|\varphi\|}_{L^2(t_\circ,t_\circ^N;V)}=1}\int_{t_\circ}^{t_\circ^N}{\langle\vartheta^\mathfrak l(t),\varphi(t)\rangle}_{V',V}\,\mathrm dt\\
&\quad\le\sup_{{\|\varphi\|}_{L^2(t_\circ,t_\circ^N;V)}=1}\int_{t_\circ}^{t_\circ^N}{\langle\rho(y(t)-y^\mathfrak l(t)),\varphi(t)\rangle}_H+{\langle\rho(y(t)^3-y^\mathfrak l(t)^3),\varphi(t)\rangle}_H\,\mathrm dt\\
&\qquad+\sup_{{\|\varphi\|}_{L^2(t_\circ,t_\circ^N;V)}=1}\int_{t_\circ}^{t_\circ^N}{\langle \mathcal P^\mathfrak l y_t(t)-y_t(t),\varphi(t)\rangle}_{V',V}-{\langle\theta \vartheta^\mathfrak l(t),\varphi(t)\rangle}_V\,\mathrm dt\\
&\quad\le \rho C_V\big({\|\varrho^\mathfrak l\|}_{L^2(t_\circ,t_\circ^N;H)}+{\|\vartheta^\mathfrak l\|}_{L^2(t_\circ,t_\circ^N;H)}\big)+\theta\,{\|\vartheta^\mathfrak l\|}_{L^2(t_\circ,t_\circ^N;V)}\\
&\qquad+C_7\,\big({\|\varrho^\mathfrak l\|}_{L^\infty(t_\circ,t_\circ^N;H)}+{\|\vartheta^\mathfrak l\|}_{L^\infty(t_\circ,t_\circ^N;H)}\big)+{\|\varrho_t^\mathfrak l\|}_{L^2(t_\circ,t_\circ^N;V')}
\end{align*}
where $C_7>0$ satisfies $C_\infty\,\|y^2+yy^\mathfrak l+(y^\mathfrak l)^2\|_{L^2(t_\circ,t_\circ^N;H)}\le C_7$. Hence, there is a constant $C_8>0$ depending on $\theta$, $\rho$, $C_W$, $C_9$, $C_{11}$ such that
\begin{equation}
\label{Eq:thetaL2Vprime}
\begin{aligned}
&{\|\vartheta^\mathfrak l\|}_{L^2(t_\circ,t_\circ^N;V')}^2\\
&\quad\le C_8\cdot\left\{
\begin{aligned}
&{\|\vartheta^\mathfrak l(t_\circ)\|}_H^2+\sum_{i=\mathfrak l+1}^d\lambda_i^V&&\text{for }X=V,\\
&{\|\vartheta^\mathfrak l(t_\circ)\|}_H^2+\sum_{i=\mathfrak l+1}^d\lambda_i^H{\|\psi_i^H-\mathcal P_H^\mathfrak l\psi_i^H\|}_V^2&&\text{for }X=H.
\end{aligned}
\right.
\end{aligned}
\end{equation}
Form \eqref{Eq:HEst}, \eqref{Eq:thetaL2V} and \eqref{Eq:thetaL2Vprime} we infer the a-priori error estimate of Theorem~\ref{Pro:Apriori}, which motivates the use of a POD approximation for our state equation \eqref{StateEquation}. \hfill$\Box$

\end{appendix}


\bibliographystyle{elsarticle-num}

\end{document}